\documentclass[reqno]{amsart}

\usepackage{amsmath,amsfonts,amssymb,amscd,verbatim,delarray,verbatim,url,epigraph06}
\usepackage[a4paper,centering,twoside,textwidth=6in]{geometry}
\usepackage[arrow,matrix]{xy}

\usepackage[russian,english]{babel}

\input cyracc.def
\newfam\cyrfam
\font\tencyr=wncyr10    scaled \magstep 2 \font\sevencyr=wncyr7
scaled \magstep 2 \font\sixcyr=wncyr6    scaled \magstep
2 2
\def\cyr{\fam\cyrfam\tencyr\cyracc}
\textfont\cyrfam=\tencyr \scriptfont\cyrfam=\sevencyr
  \scriptscriptfont\cyrfam=\sevencyr
 2 \font\tencyr=wncyr10
scaled \magstep 2  2
 2 
scaled \magstep 2  2
 2 
scaled \magstep 2

\DeclareMathOperator{\tr}{tr} \DeclareMathOperator{\Fix}{Fix}

 \DeclareMathOperator{\Sz}{Sz}
\DeclareMathOperator{\Irr}{Irr} \DeclareMathOperator{\Aut}{Aut}

\chardef\bslash=`\\ 

\begin{document}

\numberwithin{equation}{section}

\newtheorem{Theorem}{Theorem}[section]
\newtheorem{theorem}{Theorem}[section]

\newtheorem{cor}[Theorem]{Corollary}
\newtheorem{Cor}[Theorem]{Corollary}

\newtheorem{Conjecture}[Theorem]{Conjecture}
\newtheorem{conjecture}[Theorem]{Conjecture}
\newtheorem{Lemma}[Theorem]{Lemma}
\newtheorem{lemma}[Theorem]{Lemma}
\newtheorem{Property}[Theorem]{Property}
\newtheorem{property}[Theorem]{Property}
\newtheorem{condition}[Theorem]{Conditions}

\newtheorem{proposition}[Theorem]{Proposition}
\newtheorem{Proposition}[Theorem]{Proposition}
\newtheorem{Ax}[Theorem]{Axiom}
\newtheorem{Claim}[Theorem]{Claim}

\theoremstyle{definition}
\newtheorem{Definition}[Theorem]{Definition}\newtheorem{definition}[Theorem]{Definition}
\newtheorem{Problem}[Theorem]{Problem}
\newtheorem{Question}[Theorem]{Question}
\newtheorem{Example}[Theorem]{Example}
\newtheorem{example}[Theorem]{Example}
\newtheorem{remark}[Theorem]{Remark}
\newtheorem{Diagram}{Diagram}
\newtheorem{Remark}[Theorem]{Remark}

\theoremstyle{plain}

\newtheorem{conj}[Theorem]{Conjecture}
\newtheorem{corollary}[Theorem]{Corollary}
\newtheorem{prop}[Theorem]{Proposition}

\theoremstyle{definition}

\newtheorem*{quest}{Question}

\newtheorem*{answer}{Answer}
\newcommand{\diagref}[1]{diagram~\ref{#1}}
\newcommand{\thmref}[1]{Theorem~\ref{#1}}
\newcommand{\secref}[1]{Section~\ref{#1}}
\newcommand{\subsecref}[1]{Subsection~\ref{#1}}
\newcommand{\lemref}[1]{Lemma~\ref{#1}}
\newcommand{\corref}[1]{Corollary~\ref{#1}}
\newcommand{\exampref}[1]{Example~\ref{#1}}
\newcommand{\remarkref}[1]{Remark~\ref{#1}}
\newcommand{\corlref}[1]{Corollary~\ref{#1}}
\newcommand{\claimref}[1]{Claim~\ref{#1}}
\newcommand{\defnref}[1]{Definition~\ref{#1}}
\newcommand{\propref}[1]{Proposition~\ref{#1}}
\newcommand{\prref}[1]{Property~\ref{#1}}
\newcommand{\itemref}[1]{(Item~\ref{#1})}
\newcommand{\condref}[1]{Conditions~\ref{#1}}


\newcommand{\CE}{\mathcal{E}}
\newcommand{\CF}{\mathcal{F}}
\newcommand{\CG}{\mathcal{G}}\newcommand{\CV}{\mathcal{V}}
\newcommand{\CL}{\mathcal{L}}
\newcommand{\CM}{\mathcal{M}}
\newcommand{\A}{\mathcal{A}}
\newcommand{\CO}{\mathcal{O}}\newcommand{\CC}{\mathcal{C}}
\newcommand{\B}{\mathcal{B}}
\newcommand{\CS}{\mathcal{S}}
\newcommand{\CX}{\mathcal{X}}
\newcommand{\CY}{\mathcal{Y}}
\newcommand{\CT}{\mathcal{T}}
\newcommand{\CW}{\mathcal{W}}
\newcommand{\CJ}{\mathcal{J}}\newcommand{\CU}{\mathcal{U}}
\newcommand{\CR}{\mathcal{R}}

\newcommand{\st}{\sigma}
\renewcommand{\k}{\varkappa}
\newcommand{\Frac}{\mbox{Frac}}
\newcommand{\X}{\mathcal{X}}
\newcommand{\wt}{\widetilde}
\newcommand{\wh}{\widehat}
\newcommand{\mk}{\medskip}
\renewcommand{\sectionmark}[1]{}
\renewcommand{\Im}{\operatorname{Im}}
\renewcommand{\Re}{\operatorname{Re}}
\newcommand{\la}{\langle}
\newcommand{\ra}{\rangle}
\newcommand{\LND}{\mbox{LND}}
\newcommand{\Pic}{\mbox{Pic}}
\newcommand{\lnd}{\mbox{lnd}}
\newcommand{\GLND}{\mbox{GLND}}\newcommand{\glnd}{\mbox{glnd}}
\newcommand{\Der}{\mbox{DER}}\newcommand{\DER}{\mbox{DER}}
\renewcommand{\th}{\theta}
\newcommand{\ve}{\varepsilon}

\newcommand{\iy}{\infty}
\newcommand{\iintl}{\iint\limits}
\newcommand{\capl}{\operatornamewithlimits{\bigcap}\limits}
\newcommand{\cupl}{\operatornamewithlimits{\bigcup}\limits}
\newcommand{\suml}{\sum\limits}
\newcommand{\ord}{\operatorname{ord}}
\newcommand{\bk}{\bigskip}
\newcommand{\fc}{\frac}
\newcommand{\g}{\gamma}
\newcommand{\be}{\beta}
\newcommand{\dl}{\delta}
\newcommand{\Dl}{\Delta}
\newcommand{\lm}{\lambda}
\newcommand{\Lm}{\Lambda}
\newcommand{\vareps}{\varepsilon}
\newcommand{\om}{\omega}
\newcommand{\ov}{\overline}
\newcommand{\vp}{\varphi}
\newcommand{\kap}{\varkappa}
\newcommand{\al}{\alpha}
\newcommand{\Vp}{\Phi}
\newcommand{\Varphi}{\Phi}
\newcommand{\BC}{\mathbb{C}}
\newcommand{\BP}{\mathbb{P}}
\newcommand{\BQ}{\mathbb {Q}}
\newcommand{\BM}{\mathbb{M}}
\newcommand{\BR}{\mathbb{R}}\newcommand{\BN}{\mathbb{N}}
\newcommand{\BZ}{\mathbb{Z}}\newcommand{\BF}{\mathbb{F}}
\newcommand{\BA}{\mathbb {A}}

\newcommand{\BX}{\mathbb{X}}
\newcommand{\BY}{\mathbb{Y}}
\renewcommand{\Im}{\operatorname{Im}}
\newcommand{\id}{\operatorname{id}}
\newcommand{\ep}{\epsilon}
\newcommand{\tp}{\tilde\partial}
\newcommand{\doe}{\overset{\text{def}}{=}}
\newcommand{\supp} {\operatorname{supp}}
\newcommand{\loc} {\operatorname{loc}}
\newcommand{\de}{\partial}
\newcommand{\z}{\zeta}
\renewcommand{\a}{\alpha}
\newcommand{\G}{\Gamma}
\newcommand{\der}{\mbox{DER}}

\newcommand{\Spec}{\operatorname{Spec}}

\newcommand{\tG}{\widetilde G}

\newcommand{\Eq}{\Longleftrightarrow}
\newcommand{\RSL}{\operatorname{SL}}
\newcommand{\PGL}{\operatorname{PGL}}
\newcommand{\RM}{\operatorname{M}}
\newcommand{\rk}{\operatorname{rk}}
\newcommand{\codim}{\operatorname{codim}}
\newcommand{\diag}{\operatorname{diag}}
\newcommand{\ad}{\operatorname{ad}}
\newcommand{\SL}{\operatorname{SL}}
\newcommand{\GL}{\operatorname{GL}}
\newcommand{\PSL}{\operatorname{PSL}}

\newcommand{\im}{\operatorname{im}}

\newcommand{\FG}{\mathfrak {G}}
\newcommand{\FX}{\mathfrak {X}}
\newcommand{\FV}{\mathfrak {V}}
\newcommand{\Fg}{\mathfrak {g}}
\newcommand{\Fh}{\mathfrak {h}}
\newcommand{\Fb}{\mathfrak {b}}
\newcommand{\Fn}{\mathfrak {n}}
\newcommand{\Fa}{\mathfrak {a}}
\newcommand{\Fz}{\mathfrak {z}}

\newcommand{\Fsl}{\mathfrak {sl}}
\newcommand{\Fpsl}{\mathfrak {psl}}
\newcommand{\Fso}{\mathfrak {so}}
\newcommand{\Fgl}{\mathfrak {gl}}

\def\toeq{{\stackrel{\sim}{\longrightarrow}}}


\title[Equations in simple matrix groups]
{Equations in simple matrix groups: \\ algebra, geometry,
arithmetic, dynamics}

\author{Tatiana Bandman, Shelly Garion, Boris Kunyavski\u\i}

\address{Bandman: Department of
Mathematics, Bar-Ilan University, 5290002 Ramat Gan, ISRAEL}
\email{bandman@macs.biu.ac.il}

\address{Garion: Fachbereich Mathematik und Informatik, Universit\"at
M\"unster, Einsteinstrasse 62, D-48149 M\"unster, GERMANY}
\email{shelly.garion@uni-muenster.de}

\address{Kunyavski\u\i : Department of
Mathematics, Bar-Ilan University, 5290002 Ramat Gan, ISRAEL}
\email{kunyav@macs.biu.ac.il}

\begin{abstract}
We present a survey of results on word equations in simple groups,
as well as their analogues and generalizations, which were obtained
over the past decade using various methods, group-theoretic and
coming from algebraic and arithmetic geometry, number theory,
dynamical systems and computer algebra. Our focus is on
interrelations of these machineries which led to numerous
spectacular achievements, including solutions of several
long-standing problems.
\end{abstract}

\subjclass[2000]{14G05, 14G15, 20D06, 20D10, 20G40, 37P05, 37P25,
37P35, 37P55}

\keywords{matrix groups, finite simple groups, special linear group,
word map, trace map, arithmetic dynamics, periodic points, finite
fields, Lang--Weil estimate}

\maketitle

\begin{epigraph}
{{...}\\
{\cyr{\sixcyr{I snova skal{\cprime}d chuzhuyu pesnyu slozhit} \\
I kak svoyu ee proizneset.}} \\
{...}\\
{\it{Once more a skald will other's song compose \\
And utter it as of his own gift.}}} {\begin{flushright}{\cyr{\sixcyr{Osip Mandel{\cprime}shtam}}}, 1914 \\
{\it{Osip Mandelstam, 1914}} \\
(translated by Rafael Shusterovich)
\end{flushright}}
\end{epigraph}

\tableofcontents

\section{Introduction}\label{intro}

Matrix equations, which in the most general form can be written as
$$
F(A_1,\dots,A_m, X_1,\dots,X_d)=0,
$$
where $A_1,\dots,A_m$ are some fixed matrices, $X_1,\dots ,X_d$ are
unknowns, $F$ is an associative noncommutative polynomial, and the
solutions must belong to a certain class of matrices, constitute a
vast research domain, with spectacular applications well beyond
algebra, say, in areas such as differential equations and
mathematical physics. Noncommutativity brings lots of
counter-intuitive phenomena, which can be observed even when asking
deceptively simple questions about solvability of the equation or
about the number of its solutions, even within the class of
innocently looking quadratic matrix equations, see, e.g., \cite{FS},
\cite{GR1}, \cite{GR2}, \cite{EGR}, \cite{CS}, \cite{Ge}, \cite{Sl}.
(Of course, the reader coming from differential equations and
familiar, say, with Riccati matrix equations will not be too much
surprised by the difficulty of the question.)

To simplify the problem, one can limit oneself to considering
equations of the form
\begin{equation} \label{eq:matrix}
F(X_1,\dots ,X_d)=A,
\end{equation}
where $A$ is the only allowed constant matrix, and only scalar
constants are permitted to appear as coefficients of the polynomial
$F$. Even in this limited form, the question about solvability of
\eqref{eq:matrix} is far from being settled though it has been
extensively studied since the 1970's when Kaplansky asked about the
existence of polynomials whose value sets consist of the scalar
matrices. Such polynomials (called central) were discovered by
Formanek \cite{Fo} and Razmyslov \cite{Ra}; naturally, if $F$ is
central and the matrix $A$ is not scalar then equation
\eqref{eq:matrix} has no solutions. The same happens in the case
where the value set of $F$ consists of matrices with zero trace
(say, when $F$ is a Lie polynomial) and $\tr (A)\ne 0$ (and, of
course, in the trivial case where $F$ is identically zero on the
algebra  M$(n,K)$ of $n\times n$-matrices with coefficients from the
field $K$). However, there are subtler obstacles to the solvability
of \eqref{eq:matrix}, see \cite{KBMR} and references therein. Even a
very special case where $F$ is a Lie polynomial is not yet settled
though some cases where \eqref{eq:matrix} is solvable as well as
some obstacles to solvability were discussed in \cite{BGKP}.

Further simplification is essential for the present survey: we focus
our attention on the case where the solutions to the equations under
consideration must belong to a certain {\it multiplicative group} of
matrices. This naturally leads to the next modification: instead of
considering polynomials $F$ as in \eqref{eq:matrix}, which can be
viewed as elements of the free associative algebra $K\left<X_1,\dots
,X_d\right>$, we consider monomials $w(x_1,x_1^{-1},\dots,
x_d,x_d^{-1})$, which can be regarded as elements of the free
$d$-generated group $\CF_d=\left<x_1,\dots ,x_d\right>$. Thus our
main object is the following equation (where we shorten our previous
notation for $w$):
\begin{equation} \label{eq:word}
w(x_1,\dots, x_d)=g.
\end{equation}
Here $g$ is a fixed element of a group $G$, and we are looking for
solutions among $d$-tuples $(g_1,\dots ,g_d)$ of elements of $G$. By
$w(g_1,\dots ,g_d)$ we understand the element of $G$ obtained after
substituting the $g_i$ instead of the $x_i$ and performing all
multiplications and inversions in $G$. We define
$w(G)=\{w(g_1,\dots,g_d): g_1,\dots,g_d \in G\}$ as the set of
values of $w$ in $G$. Although equation \eqref{eq:word} (henceforth
called a word equation) makes sense in the case where $G$ is not
necessarily a group of matrices, we make emphasis on the matrix
case. The reason is two-fold. First, main applications we are going
to discuss refer to the cases where $G$ has nice natural matrix
representations. Second, our main goal is to stress the strength of
algebraic-geometric methods in treating problems related to word
equations. More specifically, our viewpoint can be described as
follows: let us regard all matrix entries of the $x_i$ as
indeterminates, then, after clearing denominators arising because of
determinants, we reduce equation \eqref{eq:word} to a system of
$dn^2$ polynomial equations in $n^2$ variables over the ground field
$K$ (all matrices are assumed to be of size $n\times n$). Thus we
managed to go over from a noncommutative problem to a commutative
one, but the paid price is high: the resulting system may be huge
and unaccessible to computer algebra systems. However, in certain
cases such an approach may turn out to be fruitful, and several
instances of successful implementation of this idea will follow.

Once equation \eqref{eq:word} is proclaimed as the main object of
our investigation, the first natural question one should ask oneself
is the following: what is the most natural class of groups $G$ to
start with? A strong hint is given by a theorem of Borel \cite{Bo}
(see also \cite{La}), which states that if $G$ is the group of
rational points of a connected semisimple algebraic group defined
over an infinite field, then for any $w\ne 1$ equation
\eqref{eq:word} is solvable when its right-hand side is, roughly
speaking, a ``typical'' element $g\in G$; see Section \ref{borel}
for a precise statement and some details. Thus the most substantial
part of the present survey is devoted to the case where $G$ is a
finite simple group. Although in this case Borel's theorem does not
give any direct indication which elements should be thought of as
``typical'', we shall try to convince the reader that there are
several ways out which make some sense. In particular, Larsen
\cite{La} showed that for every nontrivial word $w$ and $\ep>0$
there exists a number $C(w,\ep)$ such that if $G$ is a finite simple
group with $\#G>C(w,\ep)$ then $\#w(G) \geq \#G^{1-\ep}$, see
Section \ref{sec:image}.

One should note that there is not much hope for positive results of
Borel's flavour for groups $G$ which are too far from those singled
out in his theorem; see the discussion in \cite[Section~5]{BGKP} and
relevant references therein (to which one should add a recent paper
by Myasnikov and Nikolaev \cite{MN}). Therefore, as mentioned above,
we restrict ourselves to considering the case where $G$ runs within
a family of finite non-abelian simple groups. Moreover, since most
of the problems will be of asymptotic nature, we shall usually
ignore the sporadic groups. Being even more restrictive, we focus
our attention to the case of the family $\SL (2,q)$ (sometimes
extending it to all groups of Lie rank 1, thus adding the family of
Suzuki groups $\Sz(q)$). This case is, on the one hand, typical
enough, and sufficiently rich to have interesting applications. On
the other hand, it allows one to use some additional efficient tools
(such as trace maps) and get more precise results. Details can be
found below, in Sections \ref{WaI}, \ref{AD} and \ref{SL2}.

\medskip

One can ask various questions on equation \eqref{eq:word}. First,
one can fix a concrete group or a family of groups $G$ (say, $G=\SL
(2,q)$), a concrete element $g$ of $G$ (say, $g=1$) and a concrete
word $w$, and ask whether equation \eqref{eq:word} is solvable (or
has a non-identity solution, if $g=1$ is chosen). Even in such a
restrictive setting, some spectacular applications can be obtained,
in particular, to a long-standing problem of characterization of
finite solvable groups by recursive identities, see Section
\ref{WaI}. A variation of this approach, where one proves
solvability of a countable system of equations $w_n=1$ (with $G$
fixed as above), or, particularly, of the system $w_m=w_n$ with
$w_i$ obtained after the $i^{th}$ iteration from certain initial
data, can be reformulated in the language of periodic points. This
viewpoint was (implicitly) taken in \cite{BWW}, where another
sequence of words characterizing finite solvable groups was
constructed, and made explicit in \cite{BGKJ}. In the latter paper
this machinery has been further developed which led naturally to
some new concepts in arithmetical dynamics. Another impressive
application of arithmetical dynamics to a hard group-theoretic
problem has been demonstrated by Borisov and Sapir \cite{BS1},
\cite{BS2}, see Section \ref{Sapir}.

Further, one can consider more general questions. First, still
fixing some ``interesting'' word $w$, one may vary $g\in G$ and ask
about the solvability of \eqref{eq:word} for all $g$, or, at least,
for ``most of $g$'s''. In the language of word maps, where we denote
by
$$
w\colon G^d\to G
$$
the map given by $(g_1,\dots ,g_d)\mapsto w(g_1,\dots ,g_d)$, this
means that we are asking whether the word map is surjective or, at
least, has a ``large'' image. Ore's conjecture \cite{Or} on
representing every element of a finite simple group as a commutator
is a fabulous instance of such a setting. Due to immense work spread
over more than 50 years, it is now known that the commutator word
$w=[x,y] \in \CF_2$ satisfies $w(G)=G$ for any finite non-abelian
simple group $G$ (see \cite{LOST1}, the references therein and
Section \ref{sect.comm} below). It was therefore conjectured by
Shalev that a similar result also holds for iterated commutator
(so-called Engel words). The study of this case started in
\cite{BGG} (for $G=\SL(2,q)$, see Section \ref{sect.BG.BGG}). The
words of the form $w=x^ay^b \in \CF_2$ have also attracted special
interest. Larsen, Shalev and Tiep \cite{LST1} proved that any such
word is surjective on sufficiently large finite simple groups. By
further recent results of Guralnick and Malle \cite{GM} and of
Liebeck, O`Brien, Shalev and Tiep \cite{LOST2}, some words of the
form $x^by^b$ are known to be surjective on \emph{all} finite simple
groups. The particular case of surjectivity of these words on
$\SL(2,q)$ was studied in \cite{BG} (see Section \ref{sect.BG.BGG}).

A ramification of this sort of questions, where one asks whether any
element $g\in G$ (or ``most'' of them) can be represented as a
product of at most $k$ values of $w$ ($k$ is a fixed natural
number), has been christened ``a Waring-type problem'' by Shalev (by
analogy to a celebrated problem on representing a natural number as
a sum of $k$ powers), and a number of impressive results has been
obtained. The most conclusive one, due to Larsen, Shalev and Tiep,
says that for every nontrivial word $w$ there exists a constant
$C(w)$ such that if $G$ is a finite simple group satisfying
$\#G>C(w)$ then $w(G)^2=G$, see \cite{Sh2}, \cite{LaS2},
\cite{LST1}, \cite{LST2} and Section \ref{sect.Waring} below.

Note that a naive question whether $w(G)=G$ for any nontrivial word
$w$ and all sufficiently large finite simple non-abelian groups $G$
is clearly answered in the negative: indeed, it is easy to see that
if $G$ is a finite group and $m$ is an integer which is not
relatively prime to the order of $G$ then for the word $w=x_1^m$ one
has $w(G) \neq G$. Hence, if $v \in \CF_d$ is any word, then the
word map corresponding to $w=v^m$ cannot be surjective. A natural
question, suggested by Shalev, is whether these words are generally
the only exceptions for surjectivity of word maps in finite
non-abelian simple groups.

\medskip

Yet another type of problems, first introduced in \cite{GS}, arises
when one asks about the behaviour of the fibres of the word map
rather than about its image. Namely, Shalev asked \cite[Problem
2.10]{Sh1} whether (the cardinalities of) these fibres are
equidistributed (or close to be equidistributed) when $g$ varies in
some ``large'' subset of the image of $w$ and $G$ runs over some
family of finite groups. It was proved in \cite{GS} that the word
$w=x^2y^2 \in\CF_2$, the commutator word $w=[x,y] \in\CF_2$, as well
as the words $w =[x_1,\dots,x_d] \in \CF_d$, $d$-fold commutators in
any arrangement of brackets, are almost equidistributed on the
family of finite simple non-abelian groups. The case $G=\SL(2,q)$
was studied in some more detail in \cite{BGG}, \cite{BG}, \cite{BK},
see Section \ref{SL2}.

\medskip

We do not pretend that the present overview of word maps is
comprehensive. As one can see from the title, shamelessly
plagiarized from Manin's book \cite{Man}, it was conceived for
emphasizing the role of various tools lying outside group theory,
particularly, of algebraic-geometric and arithmetic-geometric
nature, in solving hard group-theoretic problems (note, however,
that the role of the ``algebra--geometry--arithmetic'' triad in
Manin's setting was somewhat different: his focus was on systematic
use of algebraic, mainly Galois-cohomological, methods in relating
arithmetic phenomena, such as counter-examples to local-global
principles, to geometric ones, such as non-rationality; our
viewpoint is different: a typical target is a group-theoretic
question whereas algebraic and arithmetic geometry provide an
efficient machinery). The reader inclined to purely algebraic
approaches to word equations is referred to an excellent exposition
given in a monograph by Segal \cite{Seg} and a more recent survey by
Nikolov \cite{Ni}. Vast literature on equations and system of
equations in free (and close to free) groups is left aside; see,
e.g., \cite{CRK}, particularly the introduction, for a comprehensive
bibliographical survey of this theory.

\section{Main tools} \label{sec:tools}

In this section we shall briefly describe some common machinery used
while treating word matrix equations, leaving details, which may
vary from one problem to another, for subsequent relevant sections
(note, in particular, that the problems discussed in Section
\ref{sec:word} require lots of other techniques).

\subsection{Crucial assumptions}

Our general setting is as follows.

\begin{enumerate}
\item\label{as.1}
We view a matrix group $G$ as an algebraic variety $G\subset
\BA^N_{x_1,\dots ,x_N}$ embedded into an affine space, where the
entries of a matrix $v\in G$ play role of its coordinates
$(x_1(v),\dots , x_N(v))=(v_1,\dots, v_N)$ (though this is not the
case if $G$ is a Suzuki group, $G$ then has a ``large'' subset
meeting this condition; see Section \ref{WaI}).

\item\label{as.2}
Group product $(v_1,v_2)\mapsto v=v_1v_2$ is a polynomial morphism
$G\times G\to G.$
\end{enumerate}

Changing $N$ if necessary, we can guarantee that the inversion will
also be a polynomial map; in most cases, however, this is not needed
because usually we work with matrices of determinant 1.

\medskip

Let $d\geq 2$ be an integer, and let $w=w(x_1,\dots, x_d)$ be a
nontrivial element of the free group  $\CF_d=\langle x_1,\dots, x_d
\rangle$ on $d$ letters $x_i.$  This means that it is a reduced word
in $x_i$ and $x_i^{-1}$ with nonzero exponents. Then, given a group
$H,$ the word $w$ defines a map $f_w\colon H^d\to H$,

$$f_w(h_1,\dots, h_d)=w(h_1,\dots, h_d).$$

When no confusion may arise, we will shorten $f_w$ to $w$.

Provided (\ref{as.1}) and (\ref{as.2}), we may interpret an
endomorphism of a group as a polynomial morphism $G\to G$ and a word
$w(x_1,\dots, x_d)$ as a morphism  $G^d\to G.$  Similarly, for a
word map $w\colon G^d\to G$ and an element $g\in G$, a solution of
the equation $w(x_1,\dots , x_d)=g$ corresponds to an affine
subvariety $S(w,g)\subset G^d.$

If $G$ is a {\it finite group} (and this is the case in most
applications considered in the present survey), we can choose the
ground field to be a {\it finite field} $\mathbb F$. If the problem
under consideration requires a study of a family of finite matrix
groups $G_q$, each defined over its own ground field $\mathbb F_q$,
it is convenient to view $G_q$ as a fibre of a group $\BZ$-scheme
$\CG.$ Then solutions to the equation  $w(x_1,\dots ,x_d)=g$ are
described by points $S(w,g)(\BF_q).$

This leads to the following local-global approach: instead of
solving equations in an infinite family of groups $G_q$, one has to
study a single $\BZ$-scheme $\CS(w,g)$ and then look at its fibres
at every $q.$

\subsection{Lang--Weil inequality}\label{LW}

Once the initial group-theoretic problem is reduced to an arithmetic
geometry problem of proving the existence of a rational point on a
variety defined over a finite field, a natural tool to use is the
Lang--Weil inequality \cite{LW}. It says that if $X$ is an
$n$-dimensional {\it absolutely irreducible} variety over $\BF_q$,
then asymptotically (i.e., for $q$ large enough), the number of its
rational points $\#X(\BF_q)$ does not differ too much from the
number $\#\BP^n(\BF_q)$ of points of the $n$-dimensional projective
space (which is $q^n+q^{n-1}+\dots +q+1$). Namely, the difference is
$O(q^{n-\frac{1}{2}})\leq C_1(X)q^{n-\frac{1}{2}}+C_2(X)q^{n-1}$. If
one can make this inequality {\it effective}, i.e., compute or at
least estimate $C_1$ and $C_2$, one can guarantee the existence of
an $\BF_q$-point on $X$ for $q$ large enough. Note that such a
computation is related to deep topological properties of $X$ and
requires some information on its Betti numbers (in suitable
cohomology), see \cite{GL2} for a nice survey of various
ramifications and generalizations of the Lang--Weil inequality. In
the one-dimensional case the classical Weil estimate gives an
expression of the remainder term via the genus of the curve.

One should not forget another difficulty related to checking
absolute irreducibility of $X$, which may be a highly nontrivial
task in concrete examples. Most of arising problems comprise
computational aspects, and usually require some advanced computer
algebra to overcome them, see Section \ref{WaI}.

\subsection{Trace map}\label{TM}

In the case where $G=\SL(2)$, apart from the general techniques
described in the previous section, we have at our disposal another
powerful tool going back to classical works of Vogt, Fricke and
Klein \cite{Vo}, \cite{Fr}, \cite{FK}, cited here from the paper
\cite{Ho} (see also \cite{Mag1}, \cite{Mag2}, \cite{Gol} for a nice
exposition of these results).

Let $\CF_d=\left<x_1, \dots, x_d\right>$ denote the free group on
$d$ generators.
For $G=\SL(2,k)$ ($k$ is any commutative ring with 1) and for any
$u\in \CF_d$ denote by $\tr(u)\colon G^d\to G$ the trace character,
$(g_1,\dots ,g_d)\mapsto \tr (u(g_1,\dots ,g_d))$.

\begin{Theorem}\label{Hor}
If $u$ is an arbitrary element of $\CF _d$, then the character of
$u$ can be expressed as a polynomial
$$
\tr(u)=P(t_1, \dots, t_d, t_{12},\dots ,t_{12\dots d})
$$
with integer coefficients in the $2^{d}-1$ characters
$t_{i_{1}i_{2}\dots i_{\nu}}=\tr(x_{i_{1}}x_{i_{2}}\dots
x_{i_{\nu}})$, $1\le\nu\le d, $ $1\le i_{1}<i_{2}< \dots <i_{\nu}\le
d.$
\end{Theorem}

Let $G=\SL(2,q)$, and let $\pi\colon G^d\to \BA^{2^{d}-1}$ be
defined by
$$\pi(g_1, \dots,g_d)=(t_1, \dots, t_d, t_{12},\dots ,t_{12\dots d})$$
in the notation of \thmref{Hor}.

Let $Z_d:=\pi( G^d)\subset\BA^{2^{d}-1}.$ Let $w\colon G^d\to G$ be
a word map. It follows from \thmref{Hor} that for every $d$ there
exists a polynomial map $\psi\colon \BA^{2^{d}-1}\to\BA^1$ such that
the following diagram commutes:

$$
\begin{CD}
G^d   @>{w}>>  G  \\
@V\pi VV @V\tr VV \\
Z_d(\BF_q) @>{\psi} >> \BA^1.
\end{CD}
\label{ddd}
$$

Moreover, for small $d$ we have a more precise information: one can
take $Z_2=\BA^3$ and $Z_3\subset \BA^7$ an explicitly given
hypersurface. This diagram allows one to reduce the study of the
image and fibres of $w$ to the corresponding problems for $\psi$,
which may be much simpler, see Sections \ref{tr} and \ref{3-valent}.

\section{Equations in groups of Lie rank one and characterization
of finite solvable groups} \label{WaI}

As mentioned in the introduction, even a solution of only {\it one}
word equation in a {\it small} family of finite groups may lead to
spectacular consequences. Here is an instance of such a phenomenon.

We consider a problem of characterizing various classes of groups by
identical laws. Say, $G$ is abelian if and only if $[x,y]=1$ for all
$x,y\in G$;  a finite group $G$ is nilpotent if and only if
$$
(\exists n=n(G)) \quad e_n(x,y):=
\underbrace{[[x,y],\dots,y]}_{\text{$n$ times}}\equiv 1
$$
(Zorn's theorem \cite{Zo}); what about finite {\it solvable} groups?

\subsubsection*{Two-variable sequences}
The purpose is to establish a characterization of the solvable
groups in the class of finite groups by an inductively defined,
Engel-like sequence of two-variable identities. This problem has a
long history and admits many counterparts and generalizations; the
interested reader is referred to \cite{GKP} for a survey.

\begin{Definition} \label{Sequence}
\begin{enumerate}
Let $w\in\CF_2$ and $f\in\CF_3$.
\item[(i)]
We say that a two-variable sequence $\{v_n(x,y)\}$ is defined by the
recursive law $(w,f)$ if
$$v_1(x,y)=w(x,y), \ v_{n+1}(x,y)=f(x,y,v_n(x,y)).$$

\item[(ii)]
We say that the sequence is $k$-valent  ($k=1,2,3$) if $f=f(x,y,z)$
depends on $k$ variables among which $z$ must appear.

\item[(iii)] We say that the sequence is Engel-like if $f(x,y,1)=1$ in
$\CF_3$.

\item[(iv)] We say that the sequence $\{v_n\}$ characterizes finite solvable
groups if the following holds:

A finite group $G$ is solvable if and only if
$$(\exists \  n\in\BN)  \ (\forall (x,y)\in G \times G)   \  \ v_n(x,y)=1.$$
\end{enumerate}
\end{Definition}

\begin{remark} \label{recursive}
If $\{v_n\}$ is an Engel-like sequence, the condition $v_n(x,y)=1$
implies that $v_m(x,y)=1$ for all $m\ge n$.
\end{remark}

\begin{Example}
The original Engel sequence corresponds to $w(x,y)=[x,y]$,
$f(x,y,z)=[z,y]$. It is 2-valent and Engel-like in the sense of the
previous definition but, of course, does not characterize finite
solvable groups.
\end{Example}

The first example of a sequence satisfying the conditions of
Definition \ref{Sequence} was found in \cite{BGGKPP1}. Let the
commutator be defined as $[a,b]:=aba^{-1}b^{-1}$. Set
$$w(x,y)=x^{-2}y^{-1}x, \ f(x,y,z)=[xzx^{-1},yzy^{-1}].$$ The
corresponding sequence is
\begin{equation}\label{1.1} u_1(x,y):=x^{-2}y^{-1}x, \quad
u_{n+1}(x,y):=[\, x\, u_n (x,y)\, x^{-1},\, y\, u_n(x,y)\, y^{-1}\,
].
\end{equation}

\begin{Theorem} \cite{BGGKPP1}, \cite{BGGKPP2} \label{seg1}
A finite group $G$ is solvable if and only if for some $n$ we have
$u_n(x,y)=1$ for all $x,y\in G$.
\end{Theorem}

This section is a brief exposition of the proof of this theorem. It
involves algebraic geometry, arithmetic geometry, algebra and heavy
MAGMA and SINGULAR computations.

\begin{remark}
The role of the first word is very important: for example, if one
changes the first word to a more natural one:
$u_1(x,y)=w(x,y)=[x,y]$, keeping the same recursive law $f$, we do
not know whether or not the obtained sequence characterizes finite
solvable groups. Moreover, neither do we know this for the
prototypical sequence arising from $u_1(x,y)=[x,y]$,
$f(x,y,z)=[[z,x],[z,y]]$, which characterizes finite-dimensional
solvable Lie algebras (if $[,]$ is understood as Lie bracket)
\cite{GKNP}.

More 3-valent sequences $\{u_n\}$ for which  \thmref{seg1} holds
were produced in \cite{Ri}. We conjecture after long computer
experiments that for most (if not for any) choices of $f$ satisfying
obvious necessary conditions, there is an initial word $u_1$ such
that \thmref{seg1} holds; see the discussion in Section
\ref{3-valent} and Conjecture \ref{law} below.
\end{remark}

An example of a $2$-valent sequence which still characterizes finite
solvable groups was given in \cite{BWW}:
\begin{equation}\label{eq:BWW}
s_1(x,y):=x,\quad s_{n+1}(x,y):=[ys_n(x,y)y^{-1},s_n(x,y)^{-1}].
\end{equation}

\begin{Theorem} \cite{BWW} \label{bww1}
A finite group $G$ is solvable if and only if for some $n$ we have
$s_n(x,y)=1$ for all $x,y\in G$.
\end{Theorem}

Whereas the proof of \thmref{seg1} is based on solvability of a
single word equation in a family of simple groups (see the next
section),  the proof of \thmref{bww1} requires solvability of a
system of countably many equations and has an evident dynamical
flavour, see \secref{BGKJ}.

\subsubsection*{Equations in $\SL(2)$ and Suzuki groups} 
In this section we sketch a proof of \thmref{seg1}.

Clearly in every  solvable group the identities $u_n(x,y)\equiv 1$
are satisfied from a certain $n\in\BN$ onward. By a standard
argument, using minimal counter-examples, the nontrivial direction
of \thmref{seg1} follows from

\begin{Theorem} \label{seg2}
Let $G$ be a finite non-abelian simple group. Then there are $x,\,
y\in G$ such that
\begin{equation}\label{u1u2}
u_1(x,y)=u_2(x,y), \ u_1(x,y)\ne 1.
\end{equation}
\end{Theorem}

Moreover, the same argument shows that we only need to prove
\thmref{seg2} for the groups $G$ in J.~Thompson's list of the
minimal simple non-solvable groups \cite{ThJ}. For simplicity we
slightly extend it:
\begin{enumerate}
\item\label{(1)} $G=\PSL (2,q)$ where $q\ge 4$ ($q=p^m$, $p$ prime),
\item\label{(2)} $G=\Sz (2^{2m+1})$, $m\in\BN$,
\item\label{(3)} $G=\PSL (3,3)$,
\end{enumerate}
because this does not make our proof more complicated.

Here $\Sz(2^{2m+1})$ ($m\in\BN$) denotes the Suzuki groups
(see~\cite[XI.3]{HB}).

For small groups from this list it is an easy computer exercise to
verify \thmref{seg2}. There are for example altogether $44928$
suitable pairs $x,y$ in the group $\PSL(3,3)$.

The general idea of the proof is as described in Section
\ref{sec:tools}. For a group $G$ in the above list, using a matrix
representation over $\BF _q$, we interpret solutions of the equation
$u_1(x,y)=u_2(x,y)$ as $\BF_q$-rational points of an algebraic
variety. Estimates of Lang--Weil type for the number of rational
points on a variety defined over a finite field guarantee in
appropriate circumstances the existence of such points for big $q$.
Of course we are faced here with the extra difficulty of having to
ensure that $u_1(x,y)\ne 1$ holds. This is achieved by taking $x, y$
from appropriate Zariski-open subsets only.

\subsubsection*{The case $G=\PSL(2,q)$}
Let us look for a solution to \eqref{u1u2} among the matrices of the
form
$$x(v)=\begin{bmatrix}t &-1\\1&0\end{bmatrix}, \ y(v)=\begin{bmatrix}1&b\\c&1+bc\end{bmatrix},$$
where $v=(t,b,c)\in \BA^3.$

SINGULAR computations show that the equation
$u_1(x(v),y(v))=u_2(x(v),y(v))$ defines an algebraic \emph{curve}
$C\subset \BA^3$. Note that a priori one should expect $\dim C=0$
because $C$ is defined by 3 equations in 3 variables. It is this
dimension jump that forced a somewhat peculiar choice of the initial
word $u_1(x,y)$.

It turns out that
\begin{itemize}
\item $C$ does not depend on $\BF_q;$
\item $C$ is absolutely irreducible \emph{for any $p$};
\item $p_a(\bar C)=10,  \deg  \bar C = 12$ (where $\bar C$ denotes
the projective closure of $C$, $p_a$ is the arithmetic genus and
$\deg$ is the degree).
\end{itemize}
(The second and third statements were established by computer
calculations; the proof of absolute irreducibility is technically
the most difficult part; the interested readers are referred to the
original papers.)

From Weil's estimate  for  the number of rational points over a
finite field it follows that $C(\BF_q)\ge q+1-2p_a\sqrt{q}-d$, hence
for $q>593$ there are enough rational points on $C$ to prove
\thmref{seg2}. Solutions for \eqref{u1u2} in the groups
$G=\PSL(2,q)$, $q\leq 593$, were found by computer.

\subsubsection*{The case of Suzuki groups}
In the proof of \thmref{seg2} the case of Suzuki groups $G=\Sz(q)$
($q=2^{2m+1}$) is the most difficult one.

The main reason is that although $\Sz (q)$ is contained in
$\GL(4,q)$, it is not a fibre of a $\BZ$-scheme.

In fact the group $\Sz (q)$ is defined with the help of a field
automorphism of $\BF_{q}$ (``square root of Frobenius'', see, e.g.,
\cite[Chapter 20]{Hu} for a precise definition), and hence the
standard matrix representation for $\Sz (q)$, obtained in the
original paper of Suzuki \cite{Su}, contains entries depending on
$q$. We shall describe now how our problem can still be treated by
methods of algebraic geometry.

This time one looks for a solution of \eqref{u1u2} among the
matrices parameterized by points of an $8$-dimensional space: for
$v=(a,b,c,d,a_0,b_0,c_0,d_0)\in \BA^8$, let
$$ x(v)=\begin{bmatrix}a^2a_0+ab+b_0&b&a&1\\
aa_0+b&a_0&1&0\\
a&1&0&0\\
1&0&0&0\end{bmatrix},$$
$$ y(v)=\begin{bmatrix}c^2c_0+cd+d_0&d&c&1\\
cc_0+d&c_0&1&0\\
c&1&0&0\\
1&0&0&0\end{bmatrix}.$$

Moreover, $x(v)\in \Sz(q)$ if and only if
\begin{equation}\label{xS}a_0=a^{2^{m+1}},b_0=b^{2^{m+1}}.\end{equation}
Similarly, $y(v)\in \Sz(q)$ if and only if
\begin{equation}\label{yS} c_0=c^{2^{m+1}},
d_0=d^{2^{m+1}}.\end{equation}

Define $ \Sigma:= \{v \ | \  u_1(x(v),y(v))=u_2(x(v),y(v))\}\subset
    \BA^8.$

From relations \eqref{xS} and \eqref{yS} it follows that
$$ \Sigma(\BF_{q}) =\Sigma\cap \{a_0=a^{2^m},b_0=b^{2^m},c_0=c^{2^m},
d_0=d^{2^m}\}.$$

In order to prove that the set $\Sigma(\BF_{q})$ is not empty, it is
represented as the set of fixed points of an automorphism of
$\Sigma.$

Namely, consider an automorphism $\alpha\colon \BA^8\to\BA^8$
defined by
$$\alpha (a,b,c,d,a_0,b_0,c_0,d_0)=
 (a^2_0,b_0^2,c_0^2,d_0^2,a,b,c,d).$$
Then  $\alpha^2$ is the Frobenius automorphism. Moreover, from
\eqref{xS} and \eqref{yS} it follows that $x(v),y(v)\in \Sz(q)$ if
and only if $\alpha^m(v)=v$, i.e., if and only if $v$ is a fixed
point of the $m^{th}$ iteration of $\alpha.$  Computations show that
there exists a subset $U\subset \Sigma$ such that
\begin{description}
\item[(i)] $U$ is a smooth, absolutely irreducible, affine variety,
\item[(ii)] $\dim U=2$,
\item[(iii)] $U$ is $\alpha$-invariant.
\end{description}

Denote by $b^i({U})={\rm dim} H^i_{{\text{\rm{\'et}}}}
({U},\overline{\BQ}_\ell )$ the $\ell$-adic Betti number $(\ell\neq
2)$. Then

\begin{Proposition}\label{seg3.3}
$b^1({U})\le 675$ and $b^2({U})\le 2^{22}$.
\end{Proposition}

The estimates contained in \propref{seg3.3} are derived from results
of Adolphson--Sperber \cite{AS} and Ghorpade--Lachaud \cite{GL1}
allowing one to bound the Euler characteristic of an affine variety
in terms of the number of variables, the number of defining
polynomials and their degrees. Note that since ${U}$ is affine, we
have $b^3({U})=b^4({U})=0$.  Since ${U}$ is nonsingular, the
ordinary and compact Betti numbers of  ${U}$ are related by the
Poincar\'e duality, and we have $b^i_c({U})=b^{4-i}({U})$.

In order to estimate the number $\#\Fix(U,n)$ of $\alpha^n$-fixed
points in $U,$ applied was the \emph{Lefschetz Trace Formula}:

\begin{equation} \label{LV}
\#  \Fix(U,n)=\sum_{i=0}^4(-1)^i \tr \bigl(\al^n \mid H^i_{c}
(U,{\ov{\BQ} }_{\ell } )\bigr)
\end{equation}
where $\Fix (U,n)$ is the set of fixed points of $\alpha^n$ acting
on $U$.

(A remark for the reader interested in algebraic-geometric details:
technical difficulties arising from the fact that $U$ is not
projective and $\alpha$ is not Frobenius, were overcome, roughly, as
follows:  $\alpha$ can be extended to an endomorphism of $\BP^8$
having no fixed points in $\ov U\setminus U,$ where $\ov U$ is the
projective closure of $U$. Then one can use Deligne's conjecture
stating that a formula of Lefschetz type holds after composing
$\alpha$ with a sufficiently high power of Frobenius, which, in our
case, means a high odd power of $\alpha$. Recall that Deligne's
conjecture has been proved by Fujiwara \cite{Fu}; see \cite{Va} for
simplifications and generalizations.)

From formula \eqref{LV} and Deligne's estimates for the eigenvalues
of the endomorphism induced by $\alpha$ on the \'etale cohomology,
an inequality of Lang--Weil type follows:

$$
|\#\Fix({U},n)-2^n|\leq b^1({U})\,2^{3n/4}+b^2({U})\, 2^{n/2}.
$$

An easy estimate shows that $\#{\Fix}( {U},n)\neq 0$ for $n>48$. The
cases $n<48$ were checked with the help of MAGMA.

\section{Arithmetic dynamics}\label{AD}

In this section we present two group-theoretic problems for which
the language of arithmetic dynamics appears to be an adequate one.
The first one, discussed in Section \ref{BGKJ}, arose from attempts
to understand the proof of Bray, J.~Wilson, R.~Wilson \cite{BWW},
who exhibited another sequence characterizing finite solvable groups
(see Section \ref{WaI}), find an explanation of the phenomenon that
a sequence possesses this property, and produce more such sequences.
This was essentially done in \cite{BGKJ}, where the interested
reader can find more elaborated constructions in arithmetic dynamics
which are beyond the scope of the present survey.

Another instance is related to the work of Borisov and Sapir
\cite{BS1}, \cite{BS2}, see Section \ref{Sapir}, where somewhat
similar philosophy led to an answer to another long-standing
group-theoretic question (this time, from the theory of infinite
groups).

\subsection{Verbal dynamical systems}\label{BGKJ}

Given a two-variable sequence $\{v_n(x,y)\}$ defined by a recursive
law $(w,f)$ (see Definition \ref{Sequence}), one can ask whether or
not it characterizes finite solvable groups. An obvious necessary
condition is that $\{v_n\}$ must descend along the derived series,
and we always assume that this condition holds. We also assume that
$\{v_n\}$ does not contain the identity word. Then the only
condition to check is the following one:
$$
G  {\textrm{ is not solvable }}\Rightarrow (\forall n) (\exists
(x,y)\in G\times G) : v_n(x,y)\ne 1.
$$
This condition may be further reduced (see \secref{WaI}) to the
following property.

\begin{Property}\label{t1}
Let $G$ be one of the groups $\PSL(2,q) (q\ne 2,3)$,
$\Sz(2^{2m+1})$, or $\PSL (3,3)$. Then
$$(\forall  n) (\exists (x,y)\in G\times G) : v_n(x,y)\ne 1.$$
\end{Property}

Since it is very easy to check the single case $\PSL (3,3)$, we
assume throughout below that the property holds for this group.

We will say that a sequence is \emph{very  good} if \prref{t1} holds
for all groups listed therein and is \emph{good} if it holds at
least for all $\PSL(2,q) (q\ne 2,3).$ In this section we want to
describe good sequences.

For a further simplification, the following observation is crucial:
if $\{v_n\}$ is an Engel-like sequence then Property \ref{t1} holds
in a group $G$ as soon as
$$
(\exists n) (\exists m>n) (\exists (x,y)\in G\times G) :
v_n(x,y)=v_m(x,y)\ne 1.
$$

This simple reformulation allows one to replace the proof of
solvability of a concrete equation (say, $v_1(x,y)=v_2(x,y)\ne 1$)
in a family of groups, as was done in Section \ref{WaI}, by the
proof of existence of a \emph{preperiodic point} (different from
identity) of a certain dynamical system generated by the recursive
law $f$.

Recall that a sequence $v_n(x,y)=(w(x,y),f(x,y,z))$, generated by
the first word $w$ and recursive law $f$, is formed by the rule
$$
v_1(x,y)=w(x,y), \quad v_{n+1}(x,y)=f(x,y,v_n(x,y)).
$$

In such a situation, one can define, for any group $G$, a self-map
$G\times G\times G\to G\times G\times G$ by adding ``tautological''
variables. We arrive at the notion of \emph{verbal dynamical system}
$(\CG, \tilde\phi, V)$ consisting of the following data:

\begin{itemize}
\item a group scheme $\CG$ (in our case $\CG=\SL(2,\BZ));$
\item a morphism $\phi \colon \CG\times \CG \times \CG \to \CG$
induced by the word $f(x,y,z)$;
\item\label{n2}
an endomorphism $\tilde\phi \colon \CG\times \CG\times \CG\to
\CG\times \CG \times \CG$  defined by
$\tilde\phi(x,y,z)=(x,y,\phi(x,y,z))$;
\item a forbidden set $V$ (in our case
$V= \CG\times\CG \times\{\id\});$
\item an initial word $w\colon \CG\times \CG \to \CG.$
\end{itemize}

Conversely, given such data, we reconstruct our iterative sequence
$\{v_n(x,y)\}$.

For each field  $\BF_q$ we consider the fibre
$(\SL(2,q),\tilde\phi_q,V(\BF_q)).$ The sequence $\{v_n\}$ is good
if for every $q$ there are $\tilde\phi_q$-preperiodic points outside
$V(\BF_q)$.

Consider three possible types of $f.$

\subsubsection{$1$-valent law}\label{1-valent}
The only 1-valent law is  $f(v)=v^r$. Then $v_{n+1}=v_1^{rn}.$ This
sequence is not going down along the derived series and hence does
not meet the necessary condition for characterizing solvable groups.

\subsubsection{$\PSL(2),$ $2$-valent law}\label{2-valent}

An example of a $2$-valent law is sequence \eqref{eq:BWW}. In such a
case our general setting can be simplified. Namely, if $f=f(y,z)$
does not depend on $x$, we can restrict our verbal dynamical system
to the form $\tilde\phi\colon \CG\times \CG\to\CG\times\CG$,
$(y,z)\mapsto (y,f(y,z))$, with the forbidden set $\CG\times\{\id\}$
and initial word $w(x,y)=x$.

\subsubsection{Traces}\label{tr}
In order to further simplify the dynamical system, one can use the
trace map (see~\secref{TM}). In  the special case $d=2$,
\thmref{Hor} may be formulated as follows.

\begin{Theorem}\label{trace}
Let $G=\SL(2,\BZ).$ Define $\pi\colon  G \times G  \to  \BA^3$ by $
\pi(x,y)=(\tr(x), \tr(xy), \tr(y)).$ Then for a word map $\phi\colon
G \times G  \to   G$, there is a polynomial in three variables
$P_{\phi}(s,u,t)$ such that $\tr(\phi(x,y))=P_{\phi}(\tr(x),
\tr(xy), \tr(y)).$
\end{Theorem}

Denote $s=\tr(x), u=\tr(xy), t=\tr(y).$

Let $f_1(s,u,t)=\tr(\phi(x,y))$, $f_2(s,u,t)=\tr(\phi(x,y)y)$, and
$\psi=(f_1(s,u,t),f_2(s,u,t),t)$.

According to \thmref{trace}, we have the following factorization of
the verbal dynamical system $(\CG, \tilde\phi, V)$, i.e., the
following commutative diagram:

\begin{equation}
\begin{aligned}
&  \CG \times \CG       & {} &\overset{\tilde\phi}\longrightarrow       & {} &       \CG \times \CG              \\
&\pi\downarrow  &  {}    & {}                     & {} & \downarrow\pi \\
&   \BA^3_{s,u,t} &{}      &\overset {\psi} \longrightarrow   & {}
&\BA^3_{s,u,t} .
\end{aligned}
\label{d4}
\end{equation}

In diagram (\ref{d4}),
\begin{enumerate}
\item $\pi$ is defined over $\BZ$;

\item $\pi $ is surjective  for all $\BF_q  $
 (see, e.g., \cite{Mac} or \cite{BGKJ});

\item the set $\Sigma$ of fixed points of $\psi$,
being defined by the system
$$f_1(s,u,t)=s,\ f_2(s,u,t)=u,$$
has a positive dimension.
\end{enumerate}

Respectively, for every $q$ we have a commutative diagram:
\begin{equation}\label{d5}
\begin{aligned}
&   \SL(2,q)  \times \SL(2,q)       & {} &\overset{{\tilde\phi}_q}\longrightarrow       & {} &    \SL(2,q)   \times \SL(2,q)            & \notag \\
&\pi\downarrow  &  {}    & {}                     & {} & \downarrow\pi &\notag\\
 &   \BA^3_{s,u,t}(\BF_q)   &{}      &\overset {\psi} \longrightarrow   & {} &\BA^3_{s,u,t} (\BF_q).          &\notag\end{aligned}
\end{equation}
Assume that $\Sigma(\ov{\BF}_q)$ contains an irreducible over
$\ov{\BF}_q$ curve $C$ of genus $g$ and degree $d$ that intersects
$\pi(V)$ at $k$ points. Then by Weil's inequality
$$\# (C\setminus \pi(V))(\BF_q)\geq (q+1)- 2g\sqrt{q}-d-k,$$
and $(C\setminus \pi(V))(\BF_q)\ne\emptyset $ for $q$ big enough,
$q\ge q_0(g,d,k).$ Let $a\in (C\setminus \pi(V))(\BF_q)$  and $F_a=
\pi^{-1}(a)\subset  \SL(2,q)  \times \SL(2,q) \setminus V_q.$ Then
$F_a$ is a  $\tilde\phi_q$-invariant finite set and thus contains a
nontrivial $\tilde\phi_q$-preperiodic point.

We can now formulate a procedure for checking that a $2$-valent
sequence $\{v_n\}$ is good. Namely, assume that $\{v_n\}$ is defined
by a word (law) $f(y,z).$  The process consists of the following
steps:
\begin{description}
\item[$(1)$] compute the trace map $\psi$ of the endomorphism $\tilde\phi;$
\item[$(2)$] compute the set $\Sigma $ of fixed points of $\psi;$
\item[$(3)$] find an affine curve $S\subset \Sigma$ such that

\begin{description}
\item[(i)]
$S$  is defined  over $\BZ;$
\item[(ii)]
$S$ is irreducible over $\overline{\BQ};$
\item[(iii)]
$S$ is not contained in $\pi(V).$
\end{description}

\end{description}

If the process succeeds, then one has to check finitely many
remaining cases, namely,
\begin{itemize}
\item for finitely many primes $p_1,\dots ,p_s$ one has to check whether $C$ is irreducible over $\ov\BF_{p_i}$
(cases of bad reduction);

\item for all small fields $\BF_q$ with $q<q_0(g,d,k)$ one has to find  a nontrivial $\tilde\phi$-preperiodic
point in $\SL(2,q)  \times \SL(2,q) .$\end{itemize} All these steps
can be performed by computer.

\begin{example}\label{ex-BWW}
The process was performed for the sequence $s_n.$ The curve $\Sigma$
contains a line in $\pi(V)$ and two curves, each of genus $1$ with
$3$ punctures. Thus, each of these two curves has a point over every
field $\BF_q$, $q>9.$
\end{example}

\begin{example}\label{ex-newseq}
Using this process, it was proved that the $2$-valent sequence
$$r_1(x,y)=x, \ \ r_{n+1}(x,y)=[y^2xy^{-2},x^{-1}]$$
is good as well. Moreover, it was checked for the Suzuki groups too.
Thus it characterizes finite solvable groups.
\end{example}

\subsubsection{$3$-valent sequences}\label{3-valent}
An example of a $3$-valent sequence characterizing finite solvable
groups is sequence \eqref{1.1}. Consider any $3$-valent sequence
$v_n=(w,f).$ As explained above, it gives rise to a morphism
$\phi\colon \CG^3 \to \CG$ and an endomorphism $\tilde {\phi}\colon
\CG^3  \to \CG^3$, $\tilde {\phi}(x,y,z)=(x,y,f(x,y,z)).$

The sequence is good if for every $q$ there exists $m=m(q)$ such
that there is a solution in $\SL(2,q)  \times \SL(2,q) $ to the
equation
$$v_1(x,y)=v_m(x,y)\ne 1.$$  This means that there is
a pair $(x,y)\in \SL(2,q)  \times \SL(2,q)$ such that
$$\tilde {\phi}^{m}(x,y,w(x,y))=(x,y,w(x,y))\not \in V_q,$$
where $V_q=\SL(2,q)\times \SL(2,q)\times\{\id\}$ denotes the
forbidden set.


As in \thmref{Hor}, we express the trace of $\phi(x,y,z)$ as a
polynomial in 7 variables $a_1=\tr(x)$, $a_2=\tr(y)$, $a_3=\tr(z),$
$a_{12}=\tr(xy)$, $a_{13}=\tr(xz),$ $a_{23}=\tr(yz)$,
$a_{123}=\tr(xyz)$. These variables are dependent (see, e.g.,
\cite{Mag1} or formulas (2.3)--(2.5) in \cite{Ho}):
$$
\begin{aligned}
&a_{123}^2 - a_{123}( a_{12} a_3+ a_{13} a_2+a_{23}
a_1-a_1a_2a_3)\\
& + (a_1^2+a_2^2+a_3^2+a_{12}^2+a_{13}^2+a_{23}^2 -a_1 a_2 a_{12}-
a_1 a_3 a_{13}-a_2 a_3 a_{23}+a_{12} a_{13}a_{23}-4) =0.
\end{aligned}
$$

Let $\ov{a}= (a_1, a_2,a_3, a_{12},a_{13},a_{23},a_{123})\in \BA^7$.
Let $\pi(x,u,y)=\ov{a}\in Z$  be the trace projection. Let $Z=\pi(
\CG^3)\subset  \BA^7$. Then we have a commutative diagram
\begin{equation} \label{d6}
\begin{CD}
\tG \times \tG  \times \tG   @>{\tilde\vp}>>  \tG \times \tG\times \tG  \\
@V\pi VV @V\pi VV \\
Z @>{\psi} >> Z
\end{CD}
\end{equation}
where $\psi(\ov a)=(a_1, a_2 ,l_1(\ov a), a_{12}, l_2(\ov a),l_3(\ov
a),l_4(\ov a))$,
$$l_1=\tr(\vp(x,y,z)), \ l_2=\tr(\vp(x,y,z)x),$$
$$l_3=\tr(\vp(x,y,z)y), \  l_4=\tr(\vp(x,y,z)xy).$$

In~\cite{BGKJ} it is proved that $Z$ is an irreducible hypersurface
over any algebraically closed field and that $\pi$ is surjective for
every $\BF_q.$

The dimension of the variety $F(\vp)\subset Z$ of fixed points of
$\psi$ is at least $3$. The additional condition $z=w(x,y)$ defines
a $3$-dimensional affine subset $W(w)\subset Z.$

The sequence is good if $((F(\vp)\bigcap W(w))\setminus
V)(\BF_q)\ne\emptyset $ for all $q$ big enough. Since $\dim Z=6$,
$\dim F(\vp)=3$, $\dim W(w)=3$, one should expect that
$((F(\vp)\bigcap W(w))$ is zero-dimensional. However, it turns out
that for the sequence $u_n$ defined by \eqref{1.1} it is an
absolutely irreducible curve! Thus we can formulate a sufficient
condition for a $3$-valent sequence to be good.

\begin{Theorem}\label{s3vm}
Let $w,f$ define a sequence $\{v_n(x,y)\}.$ Let $F(\vp)$ be the
variety of fixed points of the trace map $\psi$ of the corresponding
endomorphism $\tilde \varphi$  $($see diagram $(\ref{d6}))$, and let
$W(w)$ be defined by $w.$ Let $V=\{a_2=2, a_1=a_{12}, a_3=a_{23},
a_{13}=a_{123}\}$. Assume that $F(\vp)\bigcap W(w)$ contains a
positive dimensional, absolutely irreducible $\BQ$-sub\-variety
$\Phi\not\subset V.$ Then there is $q_0$ such that for every $q>q_0$
there exists a $\tilde \phi$-preperiodic point  in
$\SL(2,q)^3\setminus V_q.$
\end{Theorem}

Once again, this theorem provides a finite process, which may be
performed by computer, determining whether a sequence is good.

It is a conceptual challenge to understand whether or not the
property of a sequence to characterize finite solvable groups is
generic in some reasonable sense. We suspect that this question can
be answered in the affirmative, in the sense that for almost every
law $f(x,y,z)$, satisfying necessary conditions, there exists a
first word $w(x,y)$ such that the resulting sequence is as required.
Of course, one has to make more precise what is meant by ``almost
every law''. See Conjecture \ref{law} below.

\subsection{Mapping tori of endomorphisms of free groups}\label{Sapir}
Below we describe in brief another spectacular application of
arithmetic dynamics to group theory, following \cite{BS1}.

Given a group $G$ with generators $x_1,\dots ,x_d$ and a relation
set $R$, and its injective endomorphism $\phi$ taking $x_i$ to a
word $w_i$ ($i=1,\dots ,d$), the mapping torus $T$ of $\phi$ is
defined as the group extension of $G$ obtained by adding a generator
$t$ subject to the relations $tx_it^{-1}=w_i$. One can ask whether
$T$ is residually finite (this property means that the intersection
of the subgroups of finite index of $T$ is trivial). We refer to
\cite{Ka} and \cite{BS1} for the history and context of this
problem.

The following theorem answers this question in affirmative in the
case where the set $R$ is empty, i.e., $G$ is a free group.

\begin{theorem}\label{sapir1}
The mapping torus of any injective endomorphism of a free group is
residually finite.
\end{theorem}

Here are the main steps of the proof.

\begin{itemize}

\item
By the definition of the mapping torus of an injective endomorphism,
it is enough to prove that for any $w\in \CF_d$ and any positive
integer $a$ there is a homomorphism $h\colon T\to H$  to a finite
group $H$ such that
\begin{equation}\label{A}
h(t^awt^{-a})=h(\phi^{(a)}(w))\ne \id .
\end{equation}

\item  Any homomorphism $h_g$ of $\CF_d$ to $\SL(2)$ is
defined by a point $g=(g_1,\dots,g_d)\in \SL(2)^d.$ Then
$$h_g(x_i)=g_i, \ h_g(w(x_1,\dots,x_d))=w(g_1,\dots,g_d).$$

\item Define $\Phi\colon \SL(2)^d\to \SL(2)^d $ by
$$\Phi(s_1,\dots,s_d)=
 (w_1(s_1,\dots,s_d),\dots,w_d(s_1,\dots,s_d)).$$
Then
\begin{equation}\label{sap1}
h_g(\phi^{(a)}(w))=w(\Phi^{(a)}(g_1,\dots,g_d)).
\end{equation}

\item Take a field  $\BF_q$ big enough and consider the endomorphism
$\Phi_q\colon \SL(2,q)^d\to \SL(2,q)^d.$ Then \eqref{sap1} is still
valid.

\item  Thus any point $g\in \SL(2,q)^d$ such that

\begin{itemize}

\item it is periodic for $\Phi_q$,

\item $\pi_w(g):=w(g)=w(g_1,\dots,g_d)\ne 1$,

\end{itemize}

will have property \eqref{A}.

\item  Let $V\subset \SL(2)^d$ be the Zariski
closure of $\Phi^{(4d)}(\SL(2)^d).$ Then

\begin{itemize}

\item $V$ is  $\Phi$-invariant,

\item $\Phi\bigm|_{V}\colon V\to V$ is dominant,

\item $\pi_w(V)\ne\{\id,-\id\}.$

\end{itemize}

\item  Let $Z=V\setminus \pi_w^{-1}(\{\id, -\id\})$ and fix $\BF_q$.
By a theorem of Hrushovski \cite[Cor.~1.2]{Hr}, there is an
extension $\BF_{q_1}$ of $\BF_q$ and a $\Phi_{q_1}$-periodic point
$g\in Z(\BF_{q_1}).$
\end{itemize}

In order to avoid Hrushovski's theorem, the authors prove its
particular case for endomorphisms of affine space. That is why they
consider $\SL(2)$ as a subscheme of the scheme M${}_2$ of all
$2\times 2$ matrices. Note, however, that using Hrushovski's
theorem, they obtain a stronger result where residual finiteness is
established for the mapping torus of any endomorphism of any
finitely generated linear group. Moreover, a further refinement of
their method in \cite{BS2}, involving, in particular, search of
periodic points of self-maps defined over $p$-adic fields, allows
one to get more precise information on the structure of the mapping
tori.

\begin{remark}
Note that in the framework of the dynamical approach presented above
it is permitted to look for periodic points defined over some
extension of the original ground field. This is a subtle but
important difference from the method described in Section \ref{BGKJ}
where such an extension is forbidden, which prevents from using
Hrushovski's theorem.
\end{remark}

\section{Word maps: image and fibres}\label{sec:word}
In this section we focus on the case of finite simple groups asking
the following questions:
\begin{enumerate}
\item How big is the image of a word map?
\item What are the sizes of the fibres of a word map?
\end{enumerate}

These questions are interrelated: one can prove that the image is
large by estimating the sizes of fibres.

We start, however, with recalling a seminal result by Borel
\cite{Bo} where a general answer to the first question was obtained
for connected semisimple algebraic groups.

\subsection{Borel's theorem}\label{borel}
\begin{Theorem}\cite{Bo}. Let $G$ be a connected semisimple
algebraic group defined over a field $K$. Let $f_w\colon G^d\to G$
be the map associated to a nontrivial element $w$ of the free group
on $d\geq 2$ letters. Then $f_w$ is dominant, i.e., its image is
Zariski dense in $G$.
\end{Theorem}

Here are the main steps of the proof.

\begin{enumerate}
\item Since dominance is preserved by field extension, one may
assume that $K$ is algebraically closed of arbitrary transcendence
degree.

\item One can show that the assertion of the theorem does not depend
on the choice of $G$ within its isogeny class, so one may assume $G$
is simply connected. One can also easily reduce to the case where
$G$ is simple.

\item There is a $d$-tuple $(g_1,\dots ,g_d)\in G^d$ such that
$f_w(g_1,\dots ,g_d)\ne 1$ (because a simple group has no identical
relations).

\item First consider the case $G=\SL_n$.
\begin{itemize}
\item The following observation, going back to the unitary trick of
Weyl and used in \cite{DS}, is crucial: one can find a subfield
$L\subset K$ and a division $L$-algebra $D$ so that the group
$G=\SL_n$ contains  its anisotropic form $H=\SL_1(D)$, the group of
elements of $D$ of reduced norm 1, as a dense subset.

\item   One  proves that
\begin{equation}\label{b1}
\text{ if }\id\ne h\in H,\text{ then } 1\text { is not an eigenvalue
of }h.
\end{equation}

\item
Let $\overline{\im f_w}=Z$ be the closure of $\im f_w.$  One should
prove that $Z=\SL(n,K).$ Since $H$ is dense in $G$,
\begin{equation}\label{b1111}
\{0\}\ne f_w(H^d) \subset Z\cap H.
\end{equation}

\item \label {bi2}
Use induction on $n$. If $n=2$ and $Z\ne G,$ then $\dim Z\leq 2$ and
$Z$ is a union of a finite number of conjugacy classes of
(non-identity) semisimple elements and the set $U$ consisting of
unipotents.  Since $\id\in Z$ and $Z$ is irreducible, it is
contained in $U$, which contradicts \eqref{b1}. Hence, for $n=2$ the
statement is valid.

\item \label {bi}  Assume that the statement is proved for  $n\leq m-1.$
Fix a maximal torus $T\subset \SL(m,K)$. Let $T'\subset T$ be the
union of subtori consisting of elements with at least one eigenvalue
$1.$ $T'$ is a hypersurface in $T.$ By the induction hypothesis,
$Z\supseteq T'.$ On the other hand, by \eqref{b1111}, there is a
$d$-tuple $(h_1,\dots ,h_d)\in H^d$ such that $f_w(h_1,\dots
,h_d)\in T\setminus T'$. Since $\SL(m,K)^d$ is an irreducible
variety, $Z$ should be irreducible as well. Therefore, if it
contains a hypersurface $T'$ and at least one point outside $T',$ it
contains $T.$ Thus, $Z=\SL(m,K)$ since the conjugates of $T$ are
dense in $\SL(m,K)$ and $\im f_w$ is invariant under conjugation.
\end{itemize}

\item \label {bbi1}
Any simple group of rank $r$ not isogenous to $\SL_n$ contains a
subgroup of the same rank $r$ which is isogenous to a direct product
of groups $\SL_{n_i}$, and the assertion of the theorem follows from
the previous step.
\end{enumerate}

\begin{remark}
See \cite{KBKP} for an alternative proof based on Amitsur's theorem
on generic division rings \cite{Am}.
\end{remark}

\subsection{The image of the word map on finite simple groups} \label{sec:image}
From now on $G$ is a finite simple group, $w(x_1,\dots,x_d) \in
\CF_d$ is a nontrivial word, and we shorten our previous notation so
that $w\colon  G^d\to G $ denotes the corresponding word map.

It turns out that one can provide an analogue of Borel's theorem.
Although dominance does not make any sense in this context, it was
shown by Larsen \cite{La} that the image of the word map is large.

\begin{Theorem}\cite{La}
For every nontrivial word $w$ and any $\ep>0$ there exists
$N=N(w,\ep)$ such that if $G$ is a finite simple group of order
greater than $N$, then
$$\#w(G) \geq \#G^{1-\ep}.$$
\end{Theorem}

The original proof of this theorem heavily relied on the techniques
of Larsen--Pink \cite{LP} for estimating the sizes of the fibres of
$w$ in the case where $G$ is of Lie type. The case $G=A_n$ always
requires separate consideration, and sporadic groups can be ignored
whenever one restricts attention to asymptotic questions.

Later on, this theorem was reproved by Larsen and Shalev \cite{LaS2}
in the case of groups of Lie type using more traditional methods of
arithmetic geometry, such as Lefschetz's trace formula (not
including, however, Suzuki and Ree groups, successfully treated in
\cite{LP}).

In the same paper \cite{LaS2}, Larsen and Shalev considered some
other variations on estimating the size of $w(G)$. Namely, if $G$ is
a group of Lie type (different from ${\textsc{A}}_r$ or
${}^2{\textsc{A}}_r$) of {\it fixed Lie rank $r$}, they obtained an
estimate of the form
$$
\#w(G) > cr^{-1} \# G.
$$
Here $c$ is a positive absolute constant and $G$ is of order greater
than $N(w)$. For the alternating groups the estimate is slightly
weaker.

As there is no hope to establish surjectivity of the map $w$ for
arbitrary words (power maps provide easy counter-examples), one can
try to say something more about $w(G)$. The following terminology,
reminiscent of classical number theory, was introduced by Shalev.

\subsubsection{Waring type properties} \cite{Sh2, LaS2, LST1, LST2,
KN} \label{sect.Waring}

Waring's problem deals with expressing every natural number as a sum
of $f(k)$ $k^{th}$ powers for some suitable function $f$.
Noncommutative analogues of this problem were investigated during
the past years, answering the following questions: Can one write any
element of a finite simple group as a product of $f(k)$ $k^{th}$
powers? (See \cite{MZ}, \cite{SW}.) Can this result be extended by
replacing the word $x^k$ with an arbitrary nontrivial word $w$? Can
these results be improved by replacing the function $f$ with a
global (small) constant? See \cite{LST1} and the references therein,
and the most recent improvements in \cite{GT2}.

Using lots of various methods and ingenious techniques, Shalev
showed in \cite{Sh2} that for every nontrivial word $w$ there exists
a constant $N(w)$ such that if $G$ is a finite simple group of order
greater than $N(w)$ then $w(G)^3 = G.$ Two alternative proofs of
this result were recently found. The first, due to Nikolov and Pyber
\cite{NP}, is based on a recent result of Gowers \cite{Gow}, and the
second, for finite simple groups of bounded Lie rank, due to
Macpherson and Tent \cite{MT}, relies on model theory.

This result was substantially improved by Larsen, Shalev
and Tiep in \cite{LaS1, LaS2, LST1}.

\begin{Theorem}\cite{LST1}\label{thm.LST}
For any nontrivial word $w$ there exists a constant $N(w)$ such that
for all finite non-abelian simple groups $G$ of order greater than
$N(w)$ we have
\[
    w(G)^2=G.
\]
\end{Theorem}

The particular case of $w=x^k$ shows that this is the best possible
Waring type result for powers.

\begin{Conjecture}[Shalev]\cite[Conjectures 2.8 and 2.9]{Sh1}\label{conj.surjall}
Let $w \ne 1$ be a word which is not a proper power of another word.
Then there exists a number $N(w)$ such that, if $G$ is either $A_n$
or a finite simple group of Lie type of rank $n$, where $n > N(w)$,
then $w(G)=G.$
\end{Conjecture}

A recent result of Kassabov and Nikolov~\cite{KN} shows that the
assumption in~\thmref{thm.LST} that $G$ is sufficiently large cannot
be removed, even if we only require that $G = w(G)^k$ for a fixed
$k$. Indeed, it is shown in \cite{KN} that for any integer $k$ there
exist a word $w$ and a finite simple group $G$, such that $w$ is not
an identity in $G$, but $G \neq w(G)^k$. This is done by
constructing for any $n>13$ a specific word $w \in \CF_2$ such that
$w(A_n)$ consists of the identity and all $3$-cycles. The result
follows since for $n > 2k + 1$ there are elements in $A_n$ which
cannot be written as a product of less than $k + 1$ $3$-cycles.

Further examples of word maps on $\SL(2,2^n)$ whose image is very
small (consisting of the identity and a single conjugacy class) have
been constructed in \cite{Le1}. In a more recent preprint \cite{Lu},
Lubotzky proved (assuming the classification of finite simple
groups) that {\it any} given subset of a finite simple group $G$
which contains the identity and is invariant under $\Aut (G)$ can
arise as the image of some word map. In \cite{Le2}, this result was
extended to some almost simple and quasisimple groups.

\medskip

The main tools in~\cite{Sh2, LaS2, LST1} involve representation
theory, algebraic geometry and probabilistic methods. For any two
nontrivial words $w_1,w_2$ the rough idea is to construct special
conjugacy classes $C_1,C_2\subset G$ satisfying:
\begin{equation}\label{eq.wC}
C_1 \subset w_1(G),\ C_2 \subset w_2(G),
\end{equation}
and
\begin{equation}\label{eq.CG}
C_1C_2 \supseteq G\setminus\{1\}.
\end{equation}

The proof of \eqref{eq.CG} relies on the following classical result
of Frobenius. Let $C_1=s_1^G$, $C_2=s_2^G$. The number of ways to
write a group element $g \in G$ as $g = x_1x_2$, where $x_i \in
C_i$, is given by
\begin{equation}\label{eq.frob}
\frac{\#C_1\#C_2}{\#G} \sum_{\chi \in \Irr(G)}
\frac{\chi(s_1)\chi(s_2)\bar\chi(g)}{\chi(1)},
\end{equation}
where $\Irr(G)$ denotes the set of irreducible complex characters of
$G$.

The case of the alternating groups was established by Larsen and
Shalev in~\cite{LaS2}. First, it was proved that for any $w \neq 1$
there exists $N(w)$ such that if $n \geq N(w)$ then the image
$w(A_n)$ contains the conjugacy classes of permutations $s_n$ with a
few cycles (at most $23$), thus implying \eqref{eq.wC}. This is
highly non-elementary, involving algebraic geometry and results from
analytic number theory (such as weak versions of the Goldbach
Conjecture). The idea is to embed groups of the form $\SL(2,p)$ and
their products into $A_n$, basing on the fact that $\SL(2,p)$ embeds
into $A_{p+1}$, and an element of order $(p-1)/2$ in $\SL(2,p)$ has
two nontrivial cycles and two fixed points in this embedding, and
then use the following property of word maps on $\SL(2,p)$.

\begin{Theorem}\cite[Theorem 4.1]{LaS2}
For every nontrivial word $w$ there exist constants $M_w$ and $m_w$
with the following property: for every prime $p > M_w$, such that
$p-1$ is divisible neither by $4$ nor by any prime $3 \leq l \leq
m_w$, $w(\SL(2,p))$ contains an element of order $(p-1)/2$.
\end{Theorem}

Inclusion \eqref{eq.CG} can now be obtained by combining
\eqref{eq.frob} with the fact that all character values $\chi(s)$ of
a permutation $s \in S_n$ can be bounded in terms of the number of
cycles alone.

\begin{Theorem}\cite[Theorem 7.2]{LaS2}
Let $s \in S_n$ be a permutation with $k$ cycles $($including
$1$-cycles$)$. Then $$|\chi(s)| \leq 2^{k-1}k!$$ for all irreducible
characters $\chi$ of $S_n$.
\end{Theorem}

Larsen and Shalev~\cite{LaS2} also treated the case of finite simple
groups of Lie type of bounded rank. Later on, Larsen, Shalev and
Tiep~\cite{LST1} completed the proof for finite simple groups of Lie
type of arbitrary rank.

For finite simple groups of Lie type, $C_1$, $C_2$ are the classes
of suitable regular semisimple elements $s_1, s_2 \in G$ lying in
maximal tori $T_1, T_2 \subset G$. The tori $T_i$ are chosen so that
if $\chi$ is an irreducible character of $G$ such that
$\chi(s_1)\chi(s_2) \ne 0$ then $\chi$ is unipotent; moreover, there
is a small (in particular, bounded) number of such unipotent
characters. These results are obtained using the machinery of
Deligne--Lusztig, see, e.g., \cite{DM}. This implies that the number
of nonzero summands in $\eqref{eq.frob}$ is small, and moreover, the
character ratios can be bounded:

\begin{Theorem}\cite[Theorem 1.2.1]{LST1}
If $G$ is a finite quasisimple classical group over $\BF_q$ and $g
\in G$ is an element of support at least $N$, then
\[
    |\chi(g)|/\chi(1) < q^{-\sqrt{N}/481}
\]
for all $1_G \neq \chi \in \Irr(G)$.
\end{Theorem}

(Here the support of $g$ is defined as the codimension of its
largest eigenspace, see \cite[Definition~4.1.1]{LST1}. Recall that a
\emph{quasisimple} group $G$ is a perfect group such that $G/Z(G)$
is simple.)

The proof of \eqref{eq.wC} is based on geometric tools, and in
particular, on the Lang--Weil estimate, that allows one to establish
a Chebotarev Density Theorem for word maps.

\begin{Theorem}\cite[Corollary 5.3.3]{LST1} \label{cheb}
For every fixed nontrivial word $w$ and fixed integer $N$, there
exists $\delta > 0$ such that for every semisimple algebraic group
$\mathbb{G}$ of dimension less than $N$ over a finite field $\BF_q$
and every maximal torus $\mathbb{T}$ of $\mathbb{G}$ defined over
$\BF_q$, we have
\[
    \#(\mathbb{T}(\BF_q) \cap w(\mathbb{G}(\BF_q))) \geq \delta\#\mathbb{T}(\BF_q).
\]
\end{Theorem}

Hence, if $q$ is sufficiently large, there exist regular semisimple
elements $s_i \in w_i(\mathbb{G}(\BF_q))$ lying in any prescribed
maximal torus $\mathbb{T}(\BF_q)$. This itself is not enough, since
the group $G$ is of unbounded Lie rank. This obstacle is treated by
embedding groups $H$ of very small rank (such as $\SL_2$) over large
extension fields into $G$ so that $s_i \in w_i(H)$ remains regular
semisimple in $G$, and lies in the required maximal torus $T_i$ of
$G$. Clearly $s_i \in w_i(G)$ so that $w_i(G)$ contains the
conjugacy class $C_i = s_i^G$, as required.

\medskip

Similar Waring type results were obtained by Larsen, Shalev and Tiep
in \cite{LST2} for quasisimple groups.

\begin{Theorem}\cite{LST2}
For a fixed nontrivial word $w$ there exists a constant $N(w)$ such
that if $G$ is a finite quasisimple group of order greater than
$N(w)$, then $w(G)^3=G.$
\end{Theorem}

For various families of finite quasisimple groups, including covers
of alternating groups, a stronger result was proved in \cite{LST2}, namely that
$w(G)^2 = G$.
This was recently finalized by Guralnick and Tiep \cite{GT2}, who
proved that for any nontrivial word $w$ there exists $N=N(w)$ such
that $w(G)^2\supseteq G\setminus Z(G)$ for all quasisimple groups
$G$ of order greater than $N$.
Note, however, that in contrast with the case of simple groups
studied in \cite{LST1}, the equality $w(G)^2 = G$ may not hold for
all large finite quasisimple groups $G$. The nontrivial central
elements of finite quasisimple groups $G$ provide the main
obstructions (see Sections \ref{sect.2.power} and \ref{sect.BG.BGG}
for examples).

Further variations on the Waring theme, discussed in \cite {Sh1} and
\cite{LaS2}, include considering products $w_1(G)\dots w_k(G)$ and
intersections $w_1(G)\cap\dots \cap w_k(G)$ where $w_i$ are {\it
distinct} words. The latter case can be fit into the same framework
by looking at the fibre product of the word maps $w_i\colon G^d\to
G$ over $G$ (a different approach was suggested in \cite{NP}). The
reader is addressed to the original papers for details.

Some counterparts of Waring type properties discussed above can be
formulated for maps of matrix algebras induced by associative
noncommutative polynomials, see \cite{KBKP} for a survey.

\begin{remark}
In \cite{LaS3}, Larsen and Shalev obtained a general estimate for
the size $\#N_w(g)$ of the fibres of word maps: for {\it all} $w\in
\CF^d$, $w\ne 1$, there exists $\varepsilon > 0$ such that for {\it
all} finite simple groups $G$ and {\it all} $g\in G$ we have
$\#N_w(g) = O(\#G^{d-\varepsilon})$, where the implicit constant
depends only on $w$. Naturally, these universal estimates are rough
in comparison with the equidistribution results because they hold
for {\it all} nontrivial words, including power words, which are far
from being equidistributed, and also for {\it all} elements in the
group, including 1, for which the fibre may be very large (as in the
case of commutators).
\end{remark}

\begin{remark}
In a different spirit, estimates for the size of the fibres of word
maps were used in \cite{Ab} and \cite{NS1}, where they yielded new
criteria for distinguishing finite nilpotent and solvable groups.
(We thank the referee for pointing out the references mentioned
above.)
\end{remark}


\subsubsection{Commutators}\cite{Or, ThR, EG4, LOST1, LOST3, GM}
\label{sect.comm}

In this and next sections we consider the image of word maps for
some special words $w$. First note that for any primitive word $w$
(this means that $w$ is a part of a free generating set of $\CF_d$),
as well as for any word of the form $w=x_1^{e_1}\dots x_d^{e_d}f$,
where the $e_i$ are coprime and $f$ belongs to the derived group
$\CF'_d$, the induced map $w\colon G^d\to G$ is surjective for an
arbitrary group $G$ (see, e.g., \cite[3.1.1]{Seg}). The commutator
word is the first nontrivial instance of the surjectivity problem.

\begin{Theorem}[Ore's Conjecture]\label{thm.ore}
If $G$ is a finite non-abelian simple group, then every element of
$G$ is a commutator.
\end{Theorem}

In other words, for the commutator word $w=[x,y] \in \CF_2$, one has
$w(G)=G$ for any finite non-abelian simple group $G$. This statement
was originally posed in 1951 and proved by Ore himself for the
alternating groups \cite{Or}. During the years, this conjecture was
proved for various families of finite simple groups (see
\cite{LOST1} and the references therein). R.~Thompson \cite{ThR}
established it for the linear groups $\PSL(n,q)$, later Ellers and
Gordeev \cite{EG4} proved the conjecture for all finite simple
groups of Lie type defined over a field with more than $8$ elements,
and recently an impressive full stop was put by Liebeck, O'Brien,
Shalev and Tiep \cite{LOST1} who completed the proof for all finite
simple groups.

The original proofs of Ore~\cite{Or} and R.~Thompson~\cite{ThR} were
obtained by explicitly finding pairs of permutations (respectively,
matrices) whose commutator corresponded to some representative in
any given conjugacy class.

In order to complete the proof of Ore's Conjecture, Liebeck,
O'Brien, Shalev and Tiep used in \cite{LOST1} the following
classical criterion dating back to Frobenius, that the number of
ways to write an element $g$ in a finite group $G$ as a commutator
is
\begin{equation}\label{eq.irr.chi}
\#G \sum_{\chi \in \Irr(G)} \frac{\chi(g)}{\chi(1)}.
\end{equation}

Roughly speaking, it was shown that if $g$ is an element with a
small centralizer, then ${\chi(g)}/{\chi(1)}$ is small for $\chi \ne
1$, and the main contribution to the character
sum~\eqref{eq.irr.chi} comes from the trivial character $\chi=1$.
Hence, this sum is positive, so elements with small centralizers are
commutators. This is based on the Deligne--Lusztig theory,
and also on the theory of dual pairs and Weil
characters of classical groups \cite{TZ}, \cite{GT1}. For elements
whose centralizers are not small, the strategy is to reduce to
groups of Lie type of lower dimension and use induction. Namely, if
a certain element has a Jordan decomposition into several Jordan
blocks, and if it is possible to express each block as a commutator
in the smaller classical group, then clearly the original element is
itself a commutator.

Computer calculations (using MAGMA) played a significant role in the
proof of~\cite{LOST1}. Since the proof uses induction, it was
necessary to establish various base cases. The conjecture was proved
directly for many of these base cases by constructing the character
table of the relevant group. For various other groups certain
elements with prescribed Jordan forms as commutators were explicitly
constructed.

Similar methods were used in the subsequent paper~\cite{LOST3}, in
which it was shown that with a few (small) exceptions, every element
of a finite quasisimple group is a commutator, and moreover, any
such element is a product of two commutators.

\medskip

Ellers and Gordeev~\cite{EG4} have proved, for the finite simple
groups of Lie type over fields with more than $8$ elements, a
stronger conjecture, known as \emph{Thompson's Conjecture}.

\begin{Conjecture}[Thompson's Conjecture]\label{conj.Thompson}
Every finite simple group $G$ has a conjugacy class $C$ such that
$C^2 = G$.
\end{Conjecture}

Observe that Thompson's conjecture immediately implies Ore's
conjecture. Indeed, if $C^2 = G$ then $1 \in C^2$ so $C^{-1} = C$
and $G=CC^{-1}$. Hence, for any $g \in G$ there exist $x \in C$ and
$h \in G$ such that $g=x^hx^{-1}=[h,x]$, as required.

Thompson's conjecture was verified for many families of finite
simple groups, including the alternating groups and the sporadic
groups, see the introduction of~\cite{EG4}, but nevertheless it is
still very much open today.

The proof of Ellers and Gordeev is based on the following
generalization of the Gauss decomposition of matrices.

\begin{Theorem}\cite{EG1}--\cite{EG3}\label{thm.EG}
Let $G$ be a Chevalley group, and let $\Gamma$ be a group generated
by $G$ and a cyclic group $\langle \sigma \rangle$ which normalizes
$G$ in $\Gamma$ and acts as a diagonal automorphism on $G$
$($perhaps trivially$)$.

Let $\gamma = \sigma g \in \Gamma\setminus Z(\Gamma)$. If $h$ is any
fixed element in the group $H$, then there is $\tau \in G$ such that
$\tau \gamma \tau^{-1} = \sigma u_1hu_2,$ where $u_1 \in U^{-}$ and
$u_2 \in U$.
\end{Theorem}

Here $H$, $U$ and $U^{-}$ are the subgroups of $G$ defined by
\[
    H = \langle h_\alpha: \alpha \in \Pi \rangle, \
    U = \langle X_\alpha: \alpha \in \Phi^+ \rangle, \
    U^{-} = \langle X_\alpha: \alpha \in \Phi^- \rangle,
\]
where $\Phi$ is the root system corresponding to $G$ and $\Pi$
denotes the simple roots of $\Phi$. Recall that the Chevalley group
$G$ is generated by the root subgroups $X_\alpha$, $\alpha \in
\Phi$.

As a consequence of Theorem \ref{thm.EG}, one easily gets a
statement in the spirit of inclusion (\ref{eq.CG}):

\begin{corollary} \label{cor:reg}
If $h_1,h_2 \in H$ are regular semisimple elements in $G$ from a
maximal split torus and $C_1$, $C_2$ are the respective conjugacy
classes, then
\[
    C_1C_2 \supseteq G \setminus Z(G).
\]
\end{corollary}

Indeed, by \cite[Proposition~1]{EG1}, for fixed $h_1, h_2$ any $u_1
\in U^{-}$ and $u_2 \in U$ can be represented as
\[
    u_1 = v_1h_1v_1^{-1}h_1^{-1} \text{ and } u_2 =
    h_2^{-1}v_2h_2v_2^{-1},
\]
for some $v_1 \in U^{-}$ and $v_2 \in U$. Thus by \thmref{thm.EG},
for any noncentral conjugacy class $C \subset G$ one can find a
representative $c \in C$ such that
\[
    c = u_1h_1h_2u_2 = (v_1h_1v_1^{-1}h_1^{-1}) h_1h_2
    (h_2^{-1}v_2h_2v_2^{-1}) = (v_1h_1v_1^{-1})(v_2h_2v_2^{-1}).
\]

Corollary \ref{cor:reg} immediately implies Ore's conjecture for any
simple group $G$ containing a regular semisimple element $h$ in a
maximal split torus, and Thompson's conjecture if this element is in
addition real, i.e., if $h$ and $h^{-1}$ are conjugate. In
\cite{EG4} a careful analysis is done to show that such desired
elements actually exist in groups of Lie type over fields with more
than $8$ elements.

\medskip

In addition, Guralnick and Malle~\cite{GM} have extended the
aforementioned result \eqref{eq.CG} from~\cite{LST1} and proved the
following variant of Thompson's conjecture.

\begin{Theorem}\cite[Theorem 1.4]{GM}\label{thm.GM}
If $G$ is a finite non-abelian simple group, then there exist
conjugacy classes $C_1$, $C_2$ in $G$ with $$G \setminus \{1\}=
C_1C_2 .$$ Moreover, aside from $G = \PSL(2,7)$ or $\PSL(2,17)$, one
can assume that each $C_i$ consists of elements of order prime to
$6$.
\end{Theorem}

Similarly to~\cite{LST1}, the proof in~\cite{GM} also relies on
estimating the character sum \eqref{eq.frob} using the
Deligne--Lusztig theory, or for some small rank groups it is
computed directly from known character tables, to show that triples
$(x_1,x_2,g)$ of elements from specified conjugacy classes $C_i$
exist in a given group $G$. The conjugacy classes $C_i$ are chosen
so that only few irreducible characters vanish simultaneously on
these classes. These triples moreover generate $G$, since the
conjugacy classes $C_i$ were chosen so that their elements are
contained in few maximal subgroups of $G$.

\medskip

After considering the commutator word, it is natural to go over to
the Engel words, defined recursively by
\begin{equation}
  e_1(x,y)= [x,y] = xyx^{-1}y^{-1},\quad
  e_n(x,y)= [e_{n-1},y],
\end{equation}
and the corresponding Engel word maps $e_n\colon G \times G
\rightarrow G$. The following conjecture is naturally raised.

\begin{conjecture}[Shalev]\label{conj.engel.shalev}
Let $n \in \BN$, then the $n^{th}$ Engel word map is surjective for
any finite simple non-abelian group $G$.
\end{conjecture}

See Section \ref{SL2} for the cases $G=\SL(2,q)$ and $\PSL(2,q)$.

\subsubsection{Two-power words} \cite{LST1, GM, LOST2, LST2}
\label{sect.2.power}

It follows from \thmref{thm.LST}, due to Larsen, Shalev and Tiep
\cite{LST1}, that any two-power word is surjective on sufficiently
large finite simple groups (see \cite[Theorem 1.1.1 and Corollary
1.1.3]{LST1}). More precisely:

\begin{Theorem}\cite{LST1}\label{thm.LST.xayb}
Let $a,b$ be two nonzero integers. Then there exists a number $N =
N(a,b)$ such that if $G$ is a finite non-abelian simple group of
order at least $N$, then any element in $G$ can be written as
$x^ay^b$ for some $x,y \in G$.
\end{Theorem}

Furthermore, by recent results of Liebeck, O`Brien, Shalev and Tiep
\cite{LOST2} and of Guralnick and Malle \cite{GM} (see
\thmref{thm.GM}), some words of the form $x^by^b$ are known to be
surjective on \emph{all} finite simple groups.

\begin{Theorem}\cite[Corollary 1.5]{GM} \label{thm.GM.xbyb}
Let $G$ be a finite non-abelian simple group and let $b$ be either a
prime power or a power of $6$. Then any element in $G$ can be
written as $x^by^b$ for some $x,y \in G$.
\end{Theorem}

Note that in general, the word $x^by^b$ is not necessarily
surjective on \emph{all} finite simple groups. Indeed, if $b$ is a
multiple of the exponent of $G$ then necessarily $x^by^b=1$ for
every $x,y \in G$.

In addition, by another recent work of Larsen, Shalev and
Tiep~\cite{LST2}, if $G$ is a finite quasisimple group, then the
word $w=x^2y^2$ is surjective on $G$. On the other hand, if $b> 2$
then the word $w=x^by^b$ is \emph{not} surjective on infinitely many
finite quasisimple groups.

\subsection{The fibres of the word map}

Recall that estimating the sizes of the fibres of the word map
appeared as an integral part of estimating its image. In this
section we address a subtler problem trying to distinguish words for
which the fibres of the corresponding word map are of the same size,
at least approximately. First note that for certain words $w$ all
fibres of the word map $w\colon G^d\to G$ are exactly of the same
size {\it for any finite group $G$}. According to recent results of
\cite{Pu}, \cite{PP}, this holds only for primitive words. Primitive
words are asymptotically very rare (exponentially negligible, in the
terminology of \cite{KS}): if we count them among all words of fixed
length, their proportion tends to 0 exponentially fast (see, e.g.,
\cite{MS}). Another viewpoint at the set of primitive elements of
$\CF_d$ is that this set is closed in the profinite topology of
$\CF_d$ \cite{PP}.

We are interested in weaker equidistribution properties which hold
for more general words.

\subsubsection{Equidistribution and measure-preservation}\cite{GS,
LOST3}\label{sec.fibre}

For a word $w=w(x_1,\dots,x_d) \in \mathcal{F}_d$, a finite group
$G$ and some $g \in G$, we denote
\[
    N_w(g) = \{(g_1,\dots,g_d)\in G^d: w(g_1,\dots,g_d)=g\}.
\]
\begin{Definition} \label{defn:equi}
A word map $w\colon G^d \rightarrow G$ is \emph{almost
equidistributed} for a family of finite groups $\mathcal{G}$ if any
group $G \in \mathcal{G}$ contains a subset $S= S_G \subseteq G$
with the following properties:
\begin{enumerate}\renewcommand{\theenumi}{\it \roman{enumi}}
\item $\#S = \#G(1-\epsilon (G))$,
\item $\#N_w(g) = (\#G)^{d-1}(1+\epsilon (G))$ uniformly for all $g \in
S$,
\end{enumerate}
where $\epsilon(G) \rightarrow 0$ whenever $\#G \rightarrow \infty$.
\end{Definition}

\begin{Theorem}\cite[Theorem 1.5]{GS}\label{thm.GS.equi}
The commutator word $w=[x,y] \in \CF_2$ is almost equidistributed
for the family of finite simple groups.
\end{Theorem}

Note that we cannot require in this theorem that $S=G$. Indeed, it
is well known (and follows from \eqref{eq.irr.chi} above) that for
$w=[x,y]$ we have
\[
N_w(1) =  k(G) \#G
\]
where $k(G)$ is the number of conjugacy classes in $G$. Since $k(G)
\rightarrow \infty$ as $\#G \rightarrow \infty$ we see that the
fibre above $g=1$ is large and does not satisfy condition {(ii)}.

Two proofs are given in~\cite{GS} for this theorem. The first is
probabilistic whereas in the second the subsets $S$ are explicitly
constructed.

Let $P=P^G$ be the commutator distribution on $G$, namely $P(g) =
N_w(g)/\#G^2$, and let $U=U^G$ be the uniform distribution on $G$
(so $U(g) = 1/\#G$ for all $g \in G$). The probabilistic proof
bounds the $L_1$-distance
$$
||P-U||_1 = \sum_{g \in G} |P(g)-U(g)|
$$
between the probability measures above. Using Frobenius
formula~\eqref{eq.irr.chi} and the Cauchy--Schwarz inequality, it is
deduced in \cite[Proposition 1.1]{GS} that
$$
||P^G-U^G||_1 \le
(\zeta^G(2)-1)^{1/2},
$$
where
\[
\zeta^G(s) = \sum_{\chi \in \Irr(G)} \chi(1)^{-s}
\]
is the so-called \emph{Witten zeta function}.

Now \thmref{thm.GS.equi} follows from results of Liebeck and
Shalev~\cite{LiS2}, who showed that if $G$ is simple, then
\[
\zeta^G(s) \rightarrow 1 \hbox{ as } \#G \rightarrow \infty \hbox{
provided } s>1.
\]
This proof also provides some estimation of the function $\epsilon$.
Namely, $\epsilon(A_n) = O(n^{-1/2})$, and if $G$ is of Lie type of
rank $r$ over a field with $q$ elements then $\epsilon(G) =
O(q^{-r/4})$.

The second, constructive proof describes the subsets $S$ explicitly.
If $G=\PSL(2,q)$ then \thmref{thm.GS.equi} follows from
\eqref{eq.irr.chi} directly using the well-known character table of
$G$. If $G \neq \PSL(2,q)$ is a group of Lie type of bounded Lie
rank then $S=S_G$ is chosen as the set of regular semisimple
elements of $G$. If $G$ is a group of Lie type with unbounded Lie
rank, then $S=S_G$ contains elements whose centralizer is not very
large. If $G=A_n$ then $S=S_G$ is chosen as the set of permutations
in $A_n$ with at most $\sqrt{n}$ fixed points. This yields better
lower bounds on the cardinality of $S$. For example, in the
constructive proof for $A_n$ we obtain $\#S \ge (1- 2/[\sqrt{n}]!)
\#A_n$, which is much better than the lower bound
$(1-O(n^{-1/2}))\#A_n$ given by the probabilistic proof.

Similarly, it was shown in~\cite[Proposition 3]{LOST3} that the
commutator word is also almost equidistributed on the family of
finite quasisimple groups.

In~\cite[Section 7]{GS} it is proved that the property to be almost
equidistributed behaves well under direct products and compositions,
implying that the words $w = [x_1,\dots,x_d] \in \CF_d$, $d$-fold
commutators in any arrangement of brackets, are almost
equidistributed within the family of finite simple groups. Similar
methods are used in \cite[Theorem 7.1]{GS} to show that the word
$w=x^2y^2$ is almost equidistributed on finite simple groups, and
Larsen and Shalev have recently obtained similar results in more
general contexts.

\medskip

By \cite[Section 3]{GS}, any ``almost equidistributed'' word map
$w\colon G^d \rightarrow G$ (see~Definition \ref{defn:equi}) is also
``almost measure preserving'' in the following sense.

\begin{Definition}\label{def.measure}
A word map $w\colon G^d \rightarrow G$ is \emph{almost measure
preserving} for a family of finite groups $\mathcal{G}$ if every
group $G \in \mathcal{G}$ satisfies the following conditions:
\begin{enumerate}\renewcommand{\theenumi}{\roman{enumi}}
\item for every subset $Y \subseteq G$ we have
$$ \#w^{-1}(Y)/\#G^d = \#Y/\#G+ o(1); $$
\item for every subset $X \subseteq G^d$ we have
$$ \#w(X)/\#G \geq \#X/\#G^d - o(1); $$
\item in particular, if $X \subseteq G^d$ and $\#X/\#G^d = 1-o(1)$,
then almost every element $g \in G$ can be written as
$g=w(g_1,\dots,g_d)$ where $(g_1,\dots,g_d) \in X$;
\end{enumerate}
here $o(1)$ denotes a function depending only on $G$ which tends to
zero as $\#G \rightarrow \infty$.
\end{Definition}
This allows one to deduce in~\cite[Corollary 1.6]{GS} that the
commutator map is almost measure preserving on the family of finite
simple groups. Since almost all pairs of elements of a finite simple
group are generating pairs (see \cite{Di}, \cite{KL}, \cite{LiS1}),
the probability that some $g\in G$ can be represented as a
commutator $g = [x, y]$, where $x,y$ generate $G$, tends to $1$ as
$\#G \rightarrow \infty$, by~\cite[Theorem 1.7]{GS}.

It seems to be not so easy to extend equidistribution results from
commutators to more general words. See, however, the next section
where this is done in the particular cases $G=\SL (2,q)$ and
$\PSL(2,q)$.

\section{Word maps on $\SL(2,q)$ and $\PSL(2,q)$} \label{SL2}

\subsection{Surjectivity for Engel words and some positive
words}\label{sect.BG.BGG}
In \cite{BGG} the particular case of Engel words in the groups
$\PSL(2,q)$ and $\SL(2,q)$ was analyzed using the trace map method,
in an attempt to prove Conjecture~\ref{conj.engel.shalev} for
$\PSL(2,q)$. The main idea was to check the surjectivity of the
trace map on $\BA^3(\BF_q)$ instead of the surjectivity of the word
map on $\SL(2,q)$.

For any $x, y\in G=\SL(2,q)$ denote $s=\tr(x)$, $t=\tr(y)$ and
$u=\tr(xy)$, and define a morphism $\pi\colon G\times G\to
\BA_{s,u,t}^3$ by $\pi(x,y):=(s,u,t)$. Then the following diagram
commutes:
$$
\begin{aligned} {}
&    G \times G       & {} &\overset{\phi^n}\longrightarrow       & {} &        G \times G              & \notag \\
&\pi\downarrow  &  {}    & {}                     & {} & \downarrow\pi &\notag\\
 &   \BA^3_{s,u,t} &{}      &\overset {\psi^n} \longrightarrow   & {} &\BA^3_{s,u,t}.         &\notag\end{aligned}
$$

In this diagram, the maps $\phi^n$ and $\psi^n$ are defined
recursively as follows.
\begin{itemize}
\item
$\phi(x,y)=({[x,y]},y)\Longrightarrow$\\
$\psi(s,u,t)=(s_{1},t,t),\ s_1=s^2+u^2+t^2-sut-2;$
\medskip

\item
$\phi^{2}(x,y) =({[[x,y],y]},y)\Longrightarrow$\\
$\psi^2(s,u,t)=\psi(s_1,t,t)=(s_2,t,t),$\\
$s_2= r(s_1,t)=s_{1}^2+2t^2-s_{1}t^2-2;$\\
$\vdots$

\item $\phi^n(x,y)=(e_n(x,y),y)\Longrightarrow$\\
$\psi^n(s,u,t)=\psi^{n-1}(s_1,u,t)=\dots=\psi(s_{n-1},t,t)=(s_n,t,t),$\\
$s_{n}= r(s_{n-1},t)=r(r(s_{n-2},t),t)=\dots=r^{(n-1)}(s_1,t)$,
\end{itemize}
where $r(s,t):=s^2+2t^2-st^2-2$.

Let $\pm 2 \ne a \in \BF_q$. It is proven in  \cite{BGG}  
that a matrix $z \in \SL(2,q)$ with
$\tr(z)=a$ can be written as $z=e_{n+1}(x,y)$ for some $x,y \in
\SL(2,q)$ if and only if there is a solution $(s_1,t) \in \BF_q^2$
of the equation $r^{(n)}(s_1,t)=a$, that is if and only if the curve
$C_{n,a}$ defined over $\BF_q$ by the equation $r^{(n)}(s_1,t)=a$
has a rational point.

In \cite{BGG} it was shown that the curve $C_{n,a}$ is absolutely
irreducible over any finite field $\BF_q$. In addition, its genus
satisfies the inequality $g(C_{n,a})\leq  2^{n}(n-1)+1$, and it has
at most $\delta =5\cdot 2^{n}$ punctures. By Weil's inequality, for
$a\ne \pm 2$ and $q>2^{2n+3}(n-1)^2$ we have
$C_{n,a}(\BF_q)\ne\emptyset$, as required.

\medskip

This implies the following results, obtained in~\cite{BGG}.

\begin{Theorem}\label{thm.surj.engel}
The $n^{th}$ Engel word map is surjective on $\SL(2,q)\setminus
\{-\id\}$ $($and hence on $\PSL(2,q))$ provided that $q\geq q_0(n)$
is sufficiently large.
\end{Theorem}

On the other hand,

\begin{Proposition}
There is an infinite family of finite fields $\BF_q$ such that if
$n\geq n_0(q)$ is large enough, then the $n^{th}$ Engel word map is
\emph{not} surjective on $\SL(2,q)\setminus \{-\id\}.$
\end{Proposition}

\begin{Proposition}
For every odd prime power $q$ there is $n_0=n_0(q)$ such that
$e_n(x,y)\ne {-\id}$ for every $n>n_0$ and every $x,y\in \SL(2,q)$.
\end{Proposition}

Indeed, $z=e_{n+1}(x,y)=-\id$  implies that there is a solution to
the equation $r^{(n)}(s_1,0)=-2$, and then there is $c\in \BF_{q^2}$
such that $c^{2^n}=-1.$

\medskip

In certain cases, $e_n$ is \emph{always} surjective on $\PSL(2,q)$.

\begin{Proposition}
$e_n$ is surjective on $\SL(2,2^e)=\PSL(2,2^e)$.
\end{Proposition}

Indeed, if $q=2^e$, take $t=0$ and then $r(s,0)=s^2$, so
$r^n(s,0)=s^{2^n}$ is an isomorphism, implying that $e_n$ is
surjective on $\PSL(2,q)$.

\begin{Proposition}
$e_{n}$ is surjective on $\PSL(2,q)$ if $\sqrt{2}\in \BF_q, \
\sqrt{-1}\not\in \BF_q.$
\end{Proposition}

\begin{Proposition}
If $n \leq 4$ then $e_n$ is surjective for all groups $\PSL(2,q).$
\end{Proposition}

The last result is a consequence of \thmref{thm.surj.engel} and
MAGMA calculations.

\medskip

In \cite{BG} there is a precise description of the positive integers
$a,b$ and prime powers $q$ for which the word map $w(x,y)=x^ay^b$ is
surjective on the group $\PSL(2,q)$ (and $\SL(2,q)$). The proof is
based on the investigation of the trace map of positive words.

The key result is that $\tr(x^ay^b)$ is a \emph{linear} polynomial
in $u$, namely:
$$
\begin{aligned}
&\tr(x^ay^b) = u \cdot f_{a,b}(s,t) + h_{a,b}(s,t), \\
& \text{ where } f_{a,b}(s,t), h_{a,b}(s,t) \in \BF_q[s,t].
\end{aligned}
$$

Thus, if neither $a$ nor $b$ is divisible by the exponent of
$\PSL(2,q)$, then any element in $\BF_q$ can be written as
$\tr(x^ay^b)$ for some $x,y \in \SL(2,q)$. This immediately implies
that in this case, any \emph{semisimple} element (namely, $z \in
\SL(2,q)$ with $\tr(z) \neq \pm 2$) can be written as $z=x^ay^b$ for
some $x,y \in \SL(2,q)$. However, when $z$ is \emph{unipotent}
(namely, $z \neq \pm \id$ and $\tr(z) = \pm 2$) one has to be more
careful, and a detailed analysis is needed. Indeed, it may happen
that neither $a$ nor $b$ is divisible by the exponent of
$\PSL(2,q)$, but nevertheless the image of the word map $w=x^ay^b$
does not contain any unipotent. For example, the word
$w=x^{42}y^{42}$ is \emph{not} surjective on $\PSL(2,7)$ and
$\PSL(2,8)$.

In addition, it was determined when $-\id$ can be written as
$x^ay^b$ for some $x,y \in \SL(2,q)$. In particular, if $q \equiv
\pm 3 \bmod 8$, then $x^4y^4 \ne -\id$ for every $x,y \in \SL(2,q)$
(the same result was obtained independently in \cite{LST2}).

\medskip

These results demonstrate, in particular, the difference between
word maps in simple and quasisimple groups (see also the previous
discussion in Sections \ref{sect.Waring}, \ref{sect.comm} and
\ref{sect.2.power}).

\subsection{Criteria for equidistribution}\label{sect.equi}

In this section we describe some results on equidistribution of
solutions of word equations of the form $w(x,y)=g$ in the family of
finite groups $\SL(2,q)$ which were obtained in \cite{BK}. A
criterion for equidistribution in terms of the trace polynomial of
$w$ is given in \thmref{equi} below. This allows one to get an
explicit description of certain classes of words possessing the
equidistribution property and show that this property is generic
within these classes.  This result can be viewed, on the one hand,
as a refinement (in the $\SL_2$-case) of equidistribution theorems
of \cite{LP} and \cite{LST1} on general words $w$ and general
Chevalley groups $G$, and, on the other hand, as a generalization of
equidistribution theorems for some particular words: \cite{GS}
(commutator words on any finite simple $G$), \cite{BGG} (Engel words
on $\SL_2$), \cite{BG} (words of the form $w=x^ay^b$ on $\SL_2$).
Acting in the spirit of \cite{GS}, we deduce a criterion for
$w\colon \SL_2\times \SL_2\to \SL_2$ to be {\it almost
measure-preserving}. It turns out that ``good'' (equidistributed,
measure-preserving) words are essentially those whose trace
polynomial cannot be represented as a composition of two other
polynomials.

Here are precise definitions and results. We will follow the
approach to equidistribution adopted in \cite{GS} (see Section
\ref{sec.fibre} above).

\begin{definition} (cf. \cite[\S 3]{GS} and Section \ref{sec.fibre}) \label{def:equifin}
Let $f\colon X\to Y$ be a map between finite non-empty sets, and let
$\vareps >0$. We say that $f$ is {\em $\vareps$-equidistributed} if
there exists $Y'\subseteq Y$ such that
\begin{enumerate}
\item[(i)] $\#Y' > \#Y(1-\vareps )$;
\item[(ii)] $|\#f^{-1}(y)-\frac{\#X}{\#Y}| < \vareps \frac{\# X}{\# Y}$ for all $y\in Y'$.
\end{enumerate}
\end{definition}

The setting is as follows. Let a family of maps of finite sets
$P_q\colon X_q\to Y_q$ be given for every $q=p^n$. Assume that for
all sufficiently large $q$ the set $Y_q$ is non-empty. For each such
$q$ take $y\in Y_q$ and denote
$$
P_y=\{x\in X_q : P_q(x)=y\}.
$$

\begin{definition} \label{def:equimor-p}
Fix a prime $p$. With the notation as above, we say that the family
$P_q\colon X_q\to Y_q$, $q=p^n$, is {\em $p$-equidistributed} if
there exist a positive integer $n_0$ and a function
$\vareps_p\colon\BN\to\BN$ tending to $0$ as $n\to\infty$ such that
for all $q=p^n$ with $n>n_0$ the set $Y_q$ contains a subset $S_q$
with the following properties:
\begin{enumerate}
\item[(i)] $\#S_q < \vareps_p(q) \, (\#Y_q)$;
\item[(ii)] $|\#P_y-\frac{\#X_q}{\#Y_q}| < \vareps_p(q) \frac{\#X_q}{\#Y_q}$ for all $y\in Y_q\setminus S_q$.
\end{enumerate}
\end{definition}

\begin{remark}
Definition \ref{def:equimor-p} means that for $q=p^n$ large enough,
the map $X_q\to Y_q$ is $\vareps_p(q)$-equidistributed, in the sense
of Definition \ref{def:equifin}.
\end{remark}

\begin{definition} \label{def:equimor}
We say that the family $P_q\colon X_q\to Y_q$ is {\em
equidistributed} if it is $p$-equidistributed for all $p$ and there
exists a function $\vareps\colon\BN\to\BN$ tending to $0$ as
$n\to\infty$ such that for every $p$ and every $q=p^n$ large enough,
we have $\vareps_p(q)\le\vareps (q)$.
\end{definition}

Let us now consider the case where $Y_q=G_q= \SL(2,q),$
$X_q=(G_q)^2$ is a direct product of its two copies, and
$P_q=P_{w,q}\colon (G_q)^2\to G_q$ is the morphism induced by some
fixed word $w\in \CF_2.$

Accordingly, we say that $w$ is equidistributed (or
$p$-equidistributed) if so is the family of maps $P_{w,q}\colon \SL
(2,q)\times \SL (2,q)\to \SL (2,q)$ (or, in other words, if so is
the morphism $\BP_w\colon \SL_{2,\BZ}\times \SL_{2,\BZ}\to
\SL_{2,\BZ}$ of group schemes over $\BZ$).

Recall some properties of polynomials.

\begin{definition} \label{perm-pol}
Let $\BF$ be a finite field. We say that $h\in \BF[x]$ is a
permutation polynomial if the set of its values $\{h(z)\}_{z\in\BF}$
coincides with $\BF$.
\end{definition}

\begin{Theorem} \cite[Theorem 7.14]{LN} \label{th:perm}
Let $q=p^n$. A polynomial $h\in \BF_q[x]$ is a permutation
polynomial of all finite extensions of $\BF_q$ if and only if
$h=ax^{p^k}+b,$ where $ a\ne 0$ and $k$ is a non-negative integer.
\end{Theorem}

The following notions are essential for our criteria.

\begin{definition}\label{def3}
Let $F$ be a field. We say that a polynomial $P\in F[x_1,\dots,
x_n]$ is  {\it $F$-composite} if there exist  $Q\in F[x_1,\dots,
x_n], \deg Q\geq 1$, and $h\in F[z], \deg h\geq 2$, such that
$P=h\circ Q.$ Otherwise, we say that $P$ is {\it $F$-noncomposite}.
\end{definition}

Note that if $E/F$ is a separable field extension, it is known
\cite[Theorem 1 and Proposition 1]{AP} that $P$ is $F$-composite if
and only if $P$ is $E$-composite. In particular, working over
perfect ground fields, we may always assume, if needed, that $F$ is
algebraically closed.

\begin{definition}\label{def4}
Let $P\in \BZ[x_1,\dots, x_n]$. Fix a prime $p$.
\begin{itemize}
\item We say that $P$ is {\it $p$-composite} if
the reduced polynomial $P_p\in \BF_p[x_1,\dots, x_n]$ is
$\BF_p$-composite. Otherwise, we say that $P$ is {\it
$p$-noncomposite}.
\item We say that a $p$-composite polynomial $P$ is {\it $p$-special} if,
in the notation of Definition $\ref{def3}$, $P_p=h\circ Q$ where
$h\in \BF_p[x]$ is a permutation polynomial of all finite extensions
of $\BF_p.$
\end{itemize}
\end{definition}

\begin{definition}\label{def5}
We say that a polynomial  $P\in \BZ[x_1,\dots, x_n]$ is {\it almost
noncomposite} if for every prime $p$ it is either $p$-noncomposite
or $p$-special. Otherwise we say that $P$ is {\it very composite}.
\end{definition}

\begin{remark} \label{r1} If a polynomial  $P\in \BZ[x_1,\dots, x_n]$
is $\BQ$-noncomposite, it is $p$-noncomposite for all but finitely
many primes $p$ \cite[2.2.1]{BDN}. If $P\in \BZ[x_1,\dots, x_n]$ is
$\BQ$-composite, it is very composite.
\end{remark}

We can now formulate the main results of \cite{BK}.

\begin{Theorem}\label{p-equi}
Let $w\in \CF_2$. The morphism $\BP_w\colon\SL_{2,\BZ}\times
\SL_{2,\BZ}\to \SL_{2,\BZ}$ is $p$-equidistributed if and only if
the trace polynomial $f_w$ is either $p$-noncomposite  or
$p$-special.
\end{Theorem}

\begin{Theorem}\label{equi}
Let $w\in \CF_2$. The morphism $\BP_w\colon\SL_{2,\BZ}\times
\SL_{2,\BZ}\to \SL_{2,\BZ}$ is equidistributed if and only if the
trace polynomial $f_w$ is almost noncomposite.
\end{Theorem}

\begin{corollary}\label{surj}
Suppose that for each $p$ and all $n$ big enough the image of
$P_{w,p^n}\colon \SL(2,p^n)\times \SL(2,p^n)\to \SL(2,p^n)$ contains
all noncentral semisimple elements of $\SL(2,p^n)$. Then $w$ is
equidistributed.
\end{corollary}

For a given word $w\in \CF_2$, let us now consider the family of
groups $\hat G_q=\PSL(2,q)$ and the corresponding word maps
$w_q\colon \hat G_q\times \hat G_q\to \hat G_q$.

\begin{Proposition} \label{SL-PSL}
If the morphism $\BP_w\colon\SL_{2,\BZ}\times \SL_{2,\BZ}\to
\SL_{2,\BZ}$ is equidistributed $($or $p$-equi\-dis\-tri\-buted$)$,
then so is the family of maps $w_q\colon \hat G_q\times \hat G_q\to
\hat G_q$.
\end{Proposition}

Here are the main ingredients of the proofs of Theorems \ref{p-equi}
and \ref{equi}:

\begin{itemize}
\item Diagram
\begin{equation} \label{diagram}
\begin{CD}
G_q\times G_q @>w>> G_q \\
@V\pi VV            @VV\tr V \\
\BA_{s,u,t}^3(\BF_q)       @>f_w>>  \BA_z^1(\BF_q)
\end{CD}
\end{equation}
where $\pi (x,y)=(\tr (x), \tr(xy), \tr (y))$.

\item An explicit Lang--Weil estimate (Ghorpade--Lachaud \cite{GL1}):
if $H\subset \BA^3_{\BF_q}$ is an {\em absolutely irreducible}
hypersurface of degree $d$, then
$$
|\#H(\BF_q)-q^2|\leq (d-1)(d-2)q^{3/2}+12(d+4)^4q
$$
or, equivalently, $\#H(\BF_q)=q^2(1+r_1)$ with
$$
|r_1|\leq q^{-1/2}[(d-1)(d-2)+12(d+4)^4q^{-1/2}].
$$
For $q >Cd^8$ this gives $|r_1| < 1/2.$

\item A generalized Stein--Lorenzini inequality \cite{Na}:
if $f_{w,p}$ is $p$-noncomposite, then the {\em spectrum} $\sigma
(f_{w,p}),$ i.e., the set of all points $z\in \BA^1_{z}(\ov{\BF}_p)$
such that the hypersurface $H_z\subset\BA^3_{s,u,t}(\ov{\BF}_p)$,
defined by the equation $f_w(s,u,t)=z$, is reducible, contains at
most $d-1$ points, where $d=\deg f_w$. The same is true for each
$\sigma_q(f_w):=\sigma (f_{w,p})\bigcap \BF_q$. Let
$z\in\BA^1_{z}(\ov{\BF}_p)\setminus  \sigma(f_{w,p}).$ Then $H_z$ is
an irreducible hypersurface and hence satisfies the
Ghorpade--Lachaud inequality.

\item Estimates for fibres of the trace map:

\begin{Lemma}\label{BG}
Let $D(s,u,t)=(s^2-4)(t^2-4)(s^2+t^2+u^2-ust-4)$, and let
$\Delta\subset\BA^3_{s,u,t}$ be defined by the equation $D=0$. Let
$H \subset\BA^3_{s,u,t}(\ov{\BF}_p)$ be a hypersurface of degree $d$
such that $H\not\subset  \Delta.$ Then for $\pi$ from diagram
$(\ref{diagram})$ we have $\#\pi^{-1}(H)(\BF_q)=\#H(\BF_q)
q^3(1+r_2),$ where $|r_2| < Cd/q.$
\end{Lemma}

\item Estimates for the size of the value set of polynomials
\cite{Wa}: if $R$ is not a permutation polynomial for $\BF_q$,
$q=p^n$, then it is not a permutation polynomial for any extension
$\BF_{q^m}$ of $\BF_q$ and omits at least $(q^m-1)/d_1$ values of
$\BF_{q^m}$, where $d_1=\deg R$.

\end{itemize}

\subsubsection{Composite trace polynomials} \label{sec:comp}
The  goal of  this section is to describe words in two variables
whose trace polynomial is composite.  We send the curious reader to
\cite{BK} for proofs.

Throughout this section $T_n(x)$ stands for the $n^{th}$ Chebyshev
polynomial, and $D_n(x)=2T_n(x/2)$ for the $n^{th}$ Dickson
polynomial. It is well known (see, e.g., \cite[(2.2)]{LMT}) that
this polynomial satisfies $D_n(x+1/x)=x^n+1/x^n$ and is completely
determined by this functional equation.

We always assume that $w(x,y)$ is written in the form
\begin{equation} \label{canon}
w=x^{a_1}y^{b_1}\dots x^{a_r}y^{b_r}
\end{equation}
and is reduced (all integers $a_i$, $b_j$ are nonzero). We call the
number $r$ the {\it complexity} of $w$.

\begin{definition} \label{tr-sim}
We say that two reduced words $w=x^{a_1}y^{b_1}\dots x^{a_r}y^{b_r}$
and $v=x^{c_1}y^{d_1}\dots x^{c_{r'}}y^{d_{r'}}$, written in form
$\eqref{canon}$, are {\em trace-similar} if $r=r'$, the array
$\{|a_i|\}$ is a rearrangement of  $\{|c_i|\},$ and the array
$\{|b_i|\}$ is a rearrangement of  $\{|d_i|\}$.
\end{definition}

Note that if reduced words $w=x^{a_1}y^{b_1}\dots x^{a_r}y^{b_r}$
and $v=x^{c_1}y^{d_1}\dots x^{c_{r'}}y^{d_{r'}}$, written in form
$\eqref{canon}$, have the same trace polynomial, then $w$ and $v$
are trace-similar \cite {Ho}.

The following propositions are valid (see \cite{BK} for the proofs).

\begin{Proposition}\label{decom}
Let $w(x,y)=x^{a_1}y^{b_1}\dots x^{a_r}y^{b_r},$ $A=\sum a_i,\quad
B=\sum b_i.$ Assume that either $A\ne 0$ or $B\ne 0.$ Assume that
the trace polynomial $f_w(s,u,t)$ is $\BC$-composite,
$f_w(s,u,t)=h(q(s,u,t)),$ where $q\in \BC[s,u,t]$ and $h\in \BC[z]$,
$\deg h\ge 2$. Then $h=D_d(z)$ for some $d\ge 2.$
\end{Proposition}

\begin{Proposition}\label{sin}
Let $w$ be a reduced word of complexity $r$ written in form
\eqref{canon}. If its trace polynomial $f_w$ is $\BC$-composite,
$f_w(s,u,t)=h(q(s,u,t))$ where $q\in \BC[s,u,t]$  and  $h(x)=\mu
x^n+ \dots$ is a polynomial in one variable of degree $n$, then
$r=nm$ and $w$ is trace-similar to $v(x,y)^n$ where $v$ is a word of
complexity $m.$
\end{Proposition}

\begin{Proposition}\label{c=n}
Let $w(x,y)=x^ay^b\dots $ be a reduced word of complexity $n$ such
that $f_w(s,u,t)=D_n(q(s,u,t))$ for some $q$. Then
$w(x,y)=(x^ay^b)^n.$
\end{Proposition}

\begin{remark} \label{sin:big-p}
The statements of Propositions \ref{decom}, \ref{sin} and \ref{c=n}
remain valid if we replace $\BC$ by (the algebraic closure of) a
sufficiently big prime field $\BF_p$, and ``composite'' by
``$p$-composite'' ($p>p_0$ depending on $w$).
\end{remark}

Here are some concrete cases where one can get more conclusive
results.

\begin{corollary}\label{cor:prime}
Let $w(x, y) = x^ay^b\dots $ be a reduced word of prime complexity
$r$. If $p>r$ and $w$ is not $p$-equidistributed, then $w=v(x,
y)^r$.
\end{corollary}

\begin{corollary}\label{cor:square}
The word $w(x,y)=x^ay^bx^cy^d$ is either equidistributed or equal to
$(x^ay^b)^2$.
\end{corollary}

All facts mentioned above allow one to describe a class of words
within which a ``generic'' word induces the map which is almost
equidistributed. More precisely, we have the following proposition.

\begin{prop} \label{prop:gen}
Let $\CR$ be the set of words $w$ of {\em prime} complexity. Then
the set $S$ of words $w\in\CR$, such that the corresponding word
morphism $\BP_w\colon \SL_{2,\BZ}\times \SL_{2,\BZ}\to \SL_{2,\BZ}$
is $p$-equi\-dis\-tri\-bu\-ted for all but finitely many primes $p$,
is exponentially generic in $\CR$.
\end{prop}

(According to the terminology of \cite{KS}, this means that the
proportion of words from $S$ among all words from $\CR$ of fixed
length tends to 1 exponentially fast as the length tends to
infinity.)

This is proved by combining the results quoted above with the
well-known fact (see, e.g., \cite{AO}) stating that the class of
words which are proper powers of other words is exponentially
negligible.

\section{Concluding remarks and open problems}
We conclude with a brief discussion of various ramifications,
analogues and generalizations of results presented in this survey,
focusing on open problems, the list of which does not pretend to be
comprehensive and reflects the authors' taste. Most of them are
borrowed from \cite{BBGKP}, \cite{GKP}, \cite{BGKP}, \cite{BGG} and
\cite{BK}.

\subsubsection*{Engel-like words and solvability properties}

In light of Theorems \ref{seg1} and \ref{bww1} and subsequent
discussions in Section \ref{3-valent}, it is natural to ask whether
or not the property of a given sequence to characterize finite
solvable groups is generic. A possible way to express it is the
following conjecture.

\begin{conjecture}\label{law}
Let $\CR_0$ $($resp. $\CR)$ denote the class of words
$f(x,y,z)\in\CF_3$ satisfying the following condition: there exists
$w_0(x,y)\in\CF_2$ $($resp. $w)$ such that the sequence
$\{v^{(0)}_n(x,y)\}$ $($resp. $\{v_n(x,y)\})$ generated by the first
word $w_0$ $($resp. $w)$ and law $f$
\begin{itemize}
\item[(i)] does not contain the identity word,
\item[(ii)] is Engel-like,
\item[(iii)] descends along the derived series
\end{itemize}
$($resp. conditions \text{\rm{(i)--(iii)}} and the additional one:

{\text{\rm{(iv)}}} for any finite group $G$ the following holds: $G$
is solvable if and only if there is $n$ such that $v_n(x,y)\equiv 1$
in $G)$.

Then the class $\CR$ is generic within $\CR_0$ $($in the sense of
\cite{KS}, as in Section $\ref{sect.equi})$.
\end{conjecture}

We hope that algebraic-geometric approaches developed in Section
\ref{3-valent} could be useful in establishing this conjecture. One
can also try the following counter-part as a testing ground: in
Conjecture \ref{law}, replace throughout ``solvable'' with nilpotent
and ``derived series'' with ``lower central series''.

\medskip

It would be interesting to investigate further, in the spirit of
\cite{GKNP} and \cite{BBGKP}, relationship between finite groups and
finite-dimensional Lie algebras from the point of view of
solvability properties. Namely, one can put forward the following
(maybe over-optimistic) conjecture.

\begin{conjecture}\label{Lie}
Let $\{v_n(x,y)\}$ be an Engel-like sequence of words in the free
Lie algebra $\CW_2$ which characterizes finite-dimensional solvable
Lie algebras defined over fields of arbitrary characteristics
$($i.e., a finite-dimensional Lie algebra $\mathfrak g$ defined over
an arbitrary field $K$ is solvable if and only if there is $n$ such
that $v_n(x,y)\equiv 0$ in $\mathfrak g)$. Then the same sequence,
regarded as a sequence of words in the free group $\CF_2$ $($i.e.,
viewing Lie bracket as commutator$)$, characterizes finite solvable
groups.
\end{conjecture}

Note that the assumption on arbitrary characteristics is essential:
in \cite{BBGKP} there have been exhibited sequences characterizing
finite-dimensional Lie algebras defined over fields of
characteristic zero but not over fields of prime characteristic (and
not characterizing finite solvable groups). The first intriguing
example is the sequence defined by
$$
v_1(x,y)=[x,y], \quad v_{n+1}(v,y)=[[v_n(x,y),x],[v_n(x,y),y]].
$$
This example resists the algebraic-geometric approach described in
Section \ref{WaI} for purely computational reasons: the arising
equations lead to varieties which are out of range of SINGULAR and
MAGMA. Perhaps experts in computer algebra, who are able to apply
more sophisticated methods, will be more lucky.

\medskip

It is desirable to use Engel-like sequence for characterizing the
solvable radical of a finite group in the same way as the original
Engel sequence is used to characterize the nilpotent radical. A
theorem of Baer \cite{Ba} says that the nilpotent radical of a
finite group $G$ coincides with the set of elements $y\in G$ with
the following property: for every $x\in G$ there is $n=n(x,y)$ such
that $e_n(x,y)=1$.

\begin{Problem}
Exhibit an Engel-like sequence $\{v_n(x,y)\}$ such that the solvable
radical of any finite group $G$ coincides with the collection of
$y\in G$ with the property: for every $x\in G$ there is $n=n(x,y)$
such that $v_n(x,y)=1$.
\end{Problem}

A recent result of J.~Wilson \cite{Wi}, stating the existence of a
two-variable countable set of words with the required property,
gives a strong evidence that such an Engel-like sequence should
exist but does not provide any candidate. Note that even the toy
problem of characterizing the solvable radical of a
finite-dimensional Lie algebra with the help of Engel-like sequences
remains open in the case where the ground field is of positive
characteristic; see \cite{BBGKP} and \cite{GKP}. We hope that
algebraic-geometric machinery in the spirit of Section \ref{WaI} may
turn out to be useful for achieving this goal.

\subsubsection*{Borel's theorem and around}

Let $G$ be a connected semisimple algebraic group defined over an
infinite field $k$.

\begin{Problem} What can be said about the ``fine structure'' of $w(G)$?
In particular, describe $w$ such that $w(G)$ contains the set of all
semisimple elements of $G$.
\end{Problem}

A similar problem was discussed in \cite{KBMR} for associative
noncommutative polynomials on associative matrix algebras, in the
spirit of Kaplansky's problem. In such a setting, the resulting maps
may well be non-dominant, and not only for obvious reasons mentioned
in the introduction. In the paper cited above there have been
described certain classes of polynomials which are not central and
whose image contains elements with nonzero trace but the induced map
is not dominant on $2\times 2$-matrices. Here are some general
questions remaining open:

\begin{Question} \cite{KBMR}
Let $P$ be an associative, noncommutative, noncentral, {\em
multilinear} polynomial in $d$ variables whose image contains
matrices with nonzero trace. Does $P$ induce a dominant map
$M(n,K)^d\to M(n,K)$? Does there exist $P$ such that this map is not
surjective?
\end{Question}

Similar problems were discussed in \cite{BGKP} for Lie polynomials
$P$ on Chevalley Lie algebras $\mathfrak g$. It was shown that the
induced map is dominant provided $P$ is not identically zero on
$\mathfrak{sl}(2)$. We do not know whether or not the latter
assumption can be dropped:

\begin{Question}
Does there exist a Lie polynomial in $d$ variables, not identically
zero on a Chevalley Lie algebra $\mathfrak g$, such that the induced
map $\mathfrak g^d\to \mathfrak g$ is not dominant?
\end{Question}

It would be interesting to understand the situation with {\it
infinite-dimensional} simple Lie algebras (as well as with
finite-dimensional algebras of Cartan type over fields of positive
characteristic). The first question, which does not seem too
complicated, is the following one:

\begin{Question}
For which Lie algebras $\mathfrak g$ of Cartan type the map
$\mathfrak g\times \mathfrak g\to\mathfrak g$, $(x,y)\mapsto [x,y]$,
is surjective?
\end{Question}

If for some $\mathfrak g$ this question is answered in the
affirmative, one can continue by looking at the Engel maps, as in
\cite{BGKP}.

\begin{remark} \label{rem:multiple}
It would be interesting to consider a more general set-up when we
have a polynomial map $P\colon L^d\rightarrow L^s$. (This means that
we consider systems of equations rather than single equations.) In
\cite{GR} some dominance results were obtained for the multiple
commutator map $P\colon L\times L^d\rightarrow L^d$ given by the
formula $P(X, X_1, \dots, X_d) = ([X, X_1], \dots, [X, X_d])$.
\end{remark}

\begin{remark} \label{rem:doubleword}
In a similar spirit, one can consider generalized word maps $w\colon
G^d\to G^s$ on simple groups. Apart from \cite{GR}, see also a
discussion of a particular case $w=(w_1,w_2)\colon G^2\to G^2$ in
\cite[Problem~1]{BGGT}.
\end{remark}

\begin{remark} \label{Cremona}
It would be interesting to find more classes of infinite simple (or
close to simple) groups admitting some analogue of Borel's statement
in the sense that the image of the word map is ``large'', at least
for a generic word. Some such classes were discussed in the
literature: infinite symmetric groups \cite{Ly} (surjectivity of the
commutator word was established by Ore \cite{Or}), groups of
automorphisms of trees \cite{Mar} and of random graphs \cite{DT}.

We would suggest looking at Cremona groups, which share many common
properties with linear algebraic groups and in which such notions as
Zariski topology and dominance can be defined (see, e.g.,
\cite{Ser}).

\begin{Problem}
Let $w\in\CF_d$ be a nontrivial word. Is the corresponding word map
$w\colon (Cr(n,k))^d\to Cr(n,k)$ dominant?
\end{Problem}

Of course, the first case to be considered is $n=2$. Recall two
recent spectacular results on $Cr(2,k)$ answering long-standing
questions on its simplicity. It turns out that this group is at the
same time simple and non-simple: it does not contain nontrivial
normal subgroups closed in the topology mentioned above \cite{Bl}
but contains lots of abstract normal subgroups \cite{CL}.
\end{remark}

\begin{remark} \label{profinite}
In light of the previous remark, one can note that analogues of
Borel's theorem may look differently at the first glance when one
considers other classes of groups. Say, a certain analogue of
Borel's theorem for profinite groups is provided by the following
deep theorem by Nikolov and Segal \cite{NS2}, \cite{NS3}: let $G$ be
a finitely generated profinite group, let $w$ be a non-commutator
word, and let $\left<w(G)\right>$ denote the corresponding verbal
subgroup (i.e., the subgroup generated by the set $\{w(g)^{\pm 1}\},
g\in G$); then $\left<w(G)\right>$ is open in $G$. The original
proof in \cite{NS2} went through the solution of the restricted
Burnside problem \cite{Ze}.  
An alternative proof in \cite{NS3} does not rely on Zelmanov's
theorem (and, as pointed out to us by the referee, implies the
statement of the restricted Burnside problem).
\end{remark}

\begin{remark}
To prevent the reader from an overoptimistic view on Borel's
theorem, one has to note that there may be a significant gap between
dominance and surjectivity. Moreover, the image of the word map may
be very large in the Zariski topology (according to Borel's theorem)
but very small in some natural topology. See \cite{ThA} where such
word maps are constructed for real compact Lie groups. Note that the
behaviour of Engel word maps on such groups is much better: they are
all surjective \cite{ET}.
\end{remark}

\subsubsection*{Word maps on finite simple groups}
Regarding the image of the word map, we can recall here Thompson's
and Shalev's Conjectures \ref{conj.Thompson}, \ref{conj.surjall} and
\ref{conj.engel.shalev} discussed in Section \ref{sect.comm}. In
particular, one can ask whether the following variant of Shalev's
Conjecture \ref{conj.surjall}, for the family of groups $\PSL(2,q)$,
holds:

\begin{conj}[Shalev]
\footnote{After the first version of the present survey had been
submitted, a counter-example to this conjecture was constructed in
\cite{JLO}.} Assume that $w=w(x,y)\in \mathcal{F}_2$ is not of the
form $v(x,y)^m$ for some $v=v(x,y) \in \mathcal{F}_2$ and $m>1$.
Then there exists a constant $q_0(w)$ such that if $q>q_0(w)$ then
$w(G)=G$ for $G=\PSL(2,q)$.
\end{conj}

In particular, this conjecture holds for Engel words and for words
of the form $x^ay^b$ (see Section \ref{sect.BG.BGG} and \cite{BGG},
\cite{BG}). One can therefore attempt to use the trace map method
described in Section \ref{TM} to find more words satisfying the
above conjecture. In particular, the following questions are raised
(see Section \ref{3-valent} and \cite[Section~8]{BGG}).

\begin{Question}\cite{BGG}
What are the words $w=w(x,y)\in \mathcal{F}_2$ for which the
corresponding trace map $\psi(s,u,t)=(f_1(s,u,t),f_2(s,u,t),t)$ has
one of the following properties:
\begin{itemize}
\item[($\star$)] the set $\{f_1(s,u,t)=a\}$ is absolutely irreducible for almost
all $q$ and for every $a\in \mathbb{F}_q$?
\item[($\star\star$)] there exists a $\psi$-invariant plane $A$ and the curves
$\{\psi \bigm|_A=a\}$ are absolutely irreducible for for almost all
$q$ and for a general $a\in \mathbb{F}_q$?
\end{itemize}
\end{Question}


It is tempting to generalize the results on the image of some word
maps of \cite{BGG}, \cite{BG}, presented in Section
\ref{sect.BG.BGG}, as well as the criteria for almost
equidistribution of \cite{BK}, presented in Section \ref{sect.equi},
in the following directions:

\begin{itemize}
\item[(i)]
extend them from words in two letters to words in $d$ letters,
$d>2$;
\item[(ii)]
keep $d=2$ but consider arbitrary finite Chevalley groups;
\item[(iii)]
combine (i) and (ii).
\end{itemize}

Whereas in case (i) one can still hope to use trace polynomials,
which exist for any $d$, to produce criteria for almost
equidistribution, cases (ii) and (iii) require some new terms for
formulating such criteria and new tools for proving them.

Regardless of getting such criteria, it would be interesting to
compare, in the general case, the properties of having large image
and being equidistributed, in the spirit of Corollary \ref{surj}. We
dare to formulate the following conjecture.

\begin{conj} \label{conj:large}
For a fixed $p$, let $G_q$ be a family of Chevalley groups of fixed
Lie type over $\BF_q$ $(q=p^n$ varies$)$. For a fixed word $w\in
\CF_d$, $d\ge 2$, let $P_q=P_{w,q}\colon (G_q)^d\to G_q$ be the
corresponding map. Suppose that

$(**)$ for all $n$ big enough the image of $P_q$ contains all
regular semisimple elements of $G_q$.

\noindent Then the family $\{P_q\}$ is almost $p$-equidistributed.
\end{conj}

It is a challenging task to describe the words $w$ satisfying
condition (**) in Conjecture \ref{conj:large} (cf. the discussion in
\cite{LST1} after Theorem~5.3.2). Certainly, words of the form
$w=v^k$, $k\ge 2$, do not satisfy this condition. We do not know any
non-power word for which (**) does not hold.

\medskip

Other interesting problems, arising in the context of measure
preservation and primitivity, were raised by Puder and Parzanchevski
\cite{PP}. Namely, in \cite{PP} it was shown that the property of a
word $w$ to be measure-preserving within the class of all finite
groups can be detected on the family of all symmetric groups $S_n$.
Are there other natural families (say, $\PGL(n,q)$) that can be used
as such detectors? They also ask whether their results on measure
preservation can be extended to the class of compact Lie groups
(with respect to the Haar measure).

\medskip

One can try yet another direction: consider equidistribution
problems for matrix {\it algebras} and for polynomials more general
than word polynomials (see Introduction). Even the case of $2\times
2$-matrices is completely open.

\subsubsection*{Miscellaneous remarks and problems}

\begin{remark}
Recently, Shalev with his collaborators extended many results of
Waring type from finite simple groups to some simple algebraic
groups over $p$-adic integers. The case of simple Chevalley groups
over rings of integers remains widely open, see \cite{Sh3} for some
relevant questions and conjectures.
\end{remark}

\begin{remark}
A new dynamical viewpoint on the image of the word map was developed
by Schul and Shalev \cite{SS}. They showed that the random walk on
any finite simple group $G$, with respect to this image as a
generating set, has mixing time 2.
\end{remark}

\begin{remark} \label{rem:rings}
One could try to extend some of the results of this survey to the
case of matrix groups or algebras over some sufficiently good ring.
One has to be careful: say, in \cite{RR} there are examples of rings
$R$ such that not every element of $\Fsl (n,R)$ is a commutator.
\end{remark}

\begin{remark} \label{rem:nonassoc}
One can ask questions of Borel's type for other classes of algebras
(beyond groups, Lie algebras and associative algebras). The
interested reader may refer to \cite{Gor} for the case of values of
commutators and associators on alternative and Jordan algebras.
\end{remark}

\begin{remark}
It would be interesting to extend (at least part of) the methods and
results described in the present survey to matrix equations with
constant matrix coefficients, as discussed in the introduction.
\end{remark}

\noindent {\it Acknowledgements.} Bandman and Kunyavski\u\i \ were
supported in part by the Minerva Foundation through the Emmy Noether
Research Institute for Mathematics. Kunyavski\u\i \ was supported in part
by the Israel Science Foundation, grant 1207/12; a part of this work was 
done when he participated in the trimester program ``Arithmetic and Geometry''  
in the Hausdorff Research Institute for Mathematics (Bonn). Garion was supported by the
SFB~878 ``Groups, Geometry and Actions''. Support of these institutions is gratefully 
appreciated. 

We thank the referee for useful remarks.


\end{document}